\documentclass[11pt]{article}%
\usepackage{amsmath}
\usepackage{amsfonts}
\usepackage{amssymb}
\usepackage{graphicx}
\usepackage{a4wide}%
\newtheorem{theorem}{Theorem}

\newtheorem{corollary}[theorem]{Corollary}

\newtheorem{example}[theorem]{Example}

\newtheorem{proposition}[theorem]{Proposition}

\begin{document}

\title{On unconstrained optimization problems solved using CDT and triality theory}
\author{C. Z\u{a}linescu\\Institute of Mathematics \textquotedblleft Octav Mayer\textquotedblright,
Iasi, Romania}
\date{}
\maketitle

\textbf{Abstract} DY Gao solely or together with some of his
collaborators applied his Canonical duality theory (CDT) for solving
a class of unconstrained optimization problems, getting the
so-called ``triality theorems". Unfortunately, the ``double-min
duality" from these results published before 2010 revealed to be
false, even if in 2003 DY Gao announced that ``certain additional
conditions" are needed for getting it. After 2010 DY Gao together
with some of his collaborators published several papers in which
they added additional conditions for getting ``double-min" and
``double-max" dualities in the triality theorems. The aim of this
paper is to treat rigorously this kind of problems and to discuss
several results concerning the \textquotedblleft triality theory"
obtained up to now.

\section{Introduction}

In the preface of the book \emph{Canonical Duality Theory. Advances in
Mechanics and Mathematics, vol 37, Springer, Cham (2017)}, edited by DY Gao, V
Latorre and N Ruan, one says:

\textquotedblleft Canonical duality theory is a breakthrough methodological
theory that can be used not only for modeling complex systems within a unified
framework, but also for solving a large class of challenging problems in
multidisciplinary fields of engineering, mathematics, and sciences. ...

This theory is composed mainly of

(1) a canonical dual transformation, which can be used to formulate perfect
dual problems without duality gap;

(2) a complementary-dual principle, which solved the open problem in finite
elasticity and provides a unified analytical solution form for general
nonconvex/nonsmooth/discrete problems;

(3) a triality theory, which can be used to identify both global and local
optimality conditions and to develop powerful algorithms for solving
challenging problems in complex systems."

\bigskip

In the period 2009--2013 we published several papers in which we showed,
mainly providing counterexamples, that practically all results by DY Gao and
his collaborators called \textquotedblleft triality theorem" and published or
submitted until 2010 are false. Moreover, in the case in which the dual
function has one variable, we showed in \cite{VoiZal:13} that the ``double-min
duality" in the \textquotedblleft triality theorem" might be true only when
the primal function has also one variable. As a result, DY Gao and C Wu in
\cite{GaoWu:11} (and \cite{GaoWu:11b}, \cite{GaoWu:12}), for a particular
class of unconstrained problems, showed that the ``double-min duality" is true
only when the number of variables of the primal and dual functions are equal;
they treat the general case in \cite{GaoWu:12a} (and \cite{GaoWu:17}).

It is our aim in this work to present rigorously this \textquotedblleft
methodological theory\textquotedblright\ for unconstrained optimization
problems in finite dimensional spaces. It is not the most general framework,
but it covers all the situations met in the examples provided in DY Gao and
his collaborators' works on unconstrained optimization problems in finite
dimensions. We also point out some drawbacks and not convincing arguments from
some of those papers.

\section{Preliminaries}

We study the following unconstrained minimization problem

\bigskip

$(P)$ $~~\min$ $f(x)$ ~s.t. $x\in\mathbb{R}^{n}$

\bigskip

\noindent where $f:=q_{0}+V\circ q$ with $q(x):=\left(  q_{1}(x),...,q_{m}%
(x)\right)  ^{T}$, $q_{i}$ $(i\in\overline{0,m})$ being quadratic functions
defined on $\mathbb{R}^{n}$, and $V\in\Gamma$, $\Gamma:=\Gamma(\mathbb{R}%
^{m})$ being the class of proper convex lower semi\-continuous (lsc for short)
functions $g:\mathbb{R}^{m}\rightarrow\overline{\mathbb{R}}:=\mathbb{R}%
\cup\{-\infty,+\infty\}$. Recall that for $g:\mathbb{R}^{m}\rightarrow
\overline{\mathbb{R}}$, $\operatorname*{dom}g:=\{y\in\mathbb{R}^{m}\mid
g(y)<\infty\}$, and $g$ is proper when $\operatorname*{dom}g\neq\emptyset$ and
$g(y)\neq-\infty$ for $y\in\mathbb{R}^{m}$. The Fenchel conjugate $g^{\ast
}:\mathbb{R}^{m}\rightarrow\overline{\mathbb{R}}$ of the proper function
$g:\mathbb{R}^{m}\rightarrow\overline{\mathbb{R}}$ is defined by
\[
g^{\ast}(\sigma):=\sup\{\left\langle y,\sigma\right\rangle -g(y)\mid
y\in\mathbb{R}^{m}\}=\sup\{\left\langle y,\sigma\right\rangle -g(y)\mid
y\in\operatorname*{dom}g\}\quad(\sigma\in\mathbb{R}^{m}),
\]
while its sub\-differential at $y\in\operatorname*{dom}g$ is
\[
\partial g(y):=\left\{  \sigma\in\mathbb{R}^{m}\mid\left\langle y^{\prime
}-y,\sigma\right\rangle \leq g(y^{\prime})-g(y)~\forall y^{\prime}%
\in\mathbb{R}^{m}\right\}  ,
\]
and $\partial g(y):=\emptyset$ if $y\notin\operatorname*{dom}g;$ clearly,
\begin{equation}
g(y)+g^{\ast}(\sigma)\geq\left\langle y,\sigma\right\rangle ~~\wedge~~\left[
\sigma\in\partial g(y)\Longleftrightarrow g(y)+g^{\ast}(\sigma)=\left\langle
y,\sigma\right\rangle \quad\forall(y,\sigma)\in\mathbb{R}^{m}\times
\mathbb{R}^{m}\right]  . \label{r-fen}%
\end{equation}
It is well known that for $g\in\Gamma$ one has $g^{\ast}\in\Gamma$, and
$\sigma\in\partial g(y)$ iff $y\in\partial g^{\ast}(\sigma);$ moreover,
$\partial g(y)\neq\emptyset$ for every $y\in\operatorname*{ri}%
(\operatorname*{dom}g)$ and $g(\overline{y})=\inf_{y\in\mathbb{R}^{m}}g(y)$
iff $0\in\partial g(\overline{y})$. Because $q_{i}$ are quadratic functions,
$q_{i}(x):=\tfrac{1}{2}\left\langle x,A_{i}x\right\rangle -\left\langle
b_{i},x\right\rangle +c_{i}$ for $x\in\mathbb{R}^{n}$ with $A_{i}%
\in\mathfrak{S}_{n}$, $b_{i}\in\mathbb{R}^{n}$ (seen as column matrices), and
$c_{i}\in\mathbb{R}$ $(i\in\overline{0,m})$, where $\mathfrak{S}_{n}$ denotes
the set of $n\times n$ real symmetric matrices; of course, $c_{0}$ can be
taken to be $0$.

Consider the so called \textquotedblleft total complementary
function\textquotedblright\ (see \cite[p.~134]{GaoWu:17}), \textquotedblleft
Gao--Strang generalized complementary function\textquotedblright\ (see
\cite[p.~42]{GaoRuaPar:12}), \textquotedblleft extended Lagrangian" (see
\cite[p.~275]{Gao:00}, \cite{Gao:03}), associated to $(P)$%
\begin{equation}
\Xi:\mathbb{R}^{n}\times\mathbb{R}^{m}\rightarrow\overline{\mathbb{R}}%
,\quad\Xi(x,\sigma)=q_{0}(x)+\left\langle q(x),\sigma\right\rangle -V^{\ast
}(\sigma)=L(x,\sigma)-V^{\ast}(\sigma), \label{r-sm}%
\end{equation}
where $L$ is the (usual) Lagrangian associated to $(q_{k})_{k\in\overline
{0,m}}$, that is $L$ is the function%
\begin{equation}
L:\mathbb{R}^{n}\times\mathbb{R}^{m}\rightarrow\mathbb{R},\quad L(x,\sigma
):=q_{0}(x)+\left\langle q(x),\sigma\right\rangle . \label{r-L}%
\end{equation}
It follows that
\begin{equation}
\Xi(x,\sigma)=\tfrac{1}{2}\left\langle x,A(\sigma)x\right\rangle -\left\langle
b(\sigma),x\right\rangle +c(\sigma)-V^{\ast}(\sigma), \label{r-smb}%
\end{equation}
where, for $\sigma_{0}:=1$ and $\sigma:=(\sigma_{1},...,\sigma_{m})^{T}%
\in\mathbb{R}^{m},$%
\[
A(\sigma):=\sum\nolimits_{k=0}^{m}\sigma_{k}A_{k},\quad b(\sigma
):=\sum\nolimits_{k=0}^{m}\sigma_{k}b_{k},\quad c(\sigma):=\sum\nolimits_{k=0}%
^{m}\sigma_{k}c_{k};
\]
clearly, $A(\cdot)$, $b(\cdot)$, $c(\cdot)$ are affine functions. Hence,
$\Xi(\cdot,\sigma)$ is quadratic for each $\sigma\in\operatorname*{dom}%
V^{\ast}$ and $\Xi(x,\cdot)$ is concave for each $x\in\mathbb{R}^{n}$. Since
$V^{\ast\ast}:=\left(  V^{\ast}\right)  ^{\ast}=V$, from the definition of the
conjugate of $V^{\ast}$ and (\ref{r-sm}) we obtain that
\begin{equation}
f(x)=\sup_{\sigma\in\operatorname*{dom}V^{\ast}}\Xi(x,\sigma)=\sup_{\sigma
\in\operatorname*{ri}(\operatorname*{dom}V^{\ast})}\Xi(x,\sigma)\quad\forall
x\in\mathbb{R}^{n}, \label{r-fb}%
\end{equation}
because for a proper convex function $g:\mathbb{R}^{m}\rightarrow
\overline{\mathbb{R}}$ one has $g^{\ast}=(g+\iota_{\operatorname*{ri}%
(\operatorname*{dom}g)})^{\ast}$ (see \cite[p.~259]{Roc:72}), where the
indicator function $\iota_{C}:Z\rightarrow\overline{\mathbb{R}}$ of the subset
$C$ of a nonempty set $Z$ is defined by $\iota_{C}(z):=0$ for $z\in C$ and
$\iota_{C}(z):=\infty$ for $z\in Z\setminus C$. Moreover,%
\begin{gather}
\nabla_{x}\Xi(x,\sigma)=A(\sigma)x-b(\sigma),\quad\nabla_{xx}^{2}\Xi
(x,\sigma)=A(\sigma),\label{r-nxXi}\\
\partial\left(  -\Xi(x,\cdot)\right)  (\sigma)=\partial V^{\ast}(\sigma)-q(x),
\label{r-dsXi}%
\end{gather}
for all $(x,\sigma)\in\mathbb{R}^{n}\times\operatorname*{dom}V^{\ast}$. Hence,
for $(x,\sigma)\in\mathbb{R}^{n}\times\operatorname*{dom}V^{\ast}$ one has%
\begin{gather}
\nabla_{x}\Xi(x,\sigma)=0\Longleftrightarrow A(\sigma)x=b(\sigma
),\label{r-nxXi0}\\
0\in\partial\left(  -\Xi(x,\cdot)\right)  (\sigma)\Longleftrightarrow
q(x)\in\partial V^{\ast}(\sigma)\Longleftrightarrow\sigma\in\partial V\left(
(q(x)\right)  . \label{r-dsXi0}%
\end{gather}

Consider the following sets in which $\sigma$ is taken from $\mathbb{R}^{n}$
if not specified otherwise:%
\begin{gather*}
Y_{0}:=\{\sigma\mid\det A(\sigma)\neq0\},\quad Y^{+}:=\{\sigma\mid
A(\sigma)\succ0\},\quad Y^{-}:=\{\sigma\mid A(\sigma)\prec0\},\\
Y_{\operatorname{col}}:=\{\sigma\mid b(\sigma)\in\operatorname{Im}%
A(\sigma)\},~Y_{\operatorname{col}}^{+}:=\{\sigma\in Y_{\operatorname{col}%
}\mid A(\sigma)\succeq0\},~Y_{\operatorname{col}}^{+}:=\{\sigma\in
Y_{\operatorname{col}}\mid A(\sigma)\preceq0\},\\
S_{0}:=Y_{0}\cap\operatorname*{dom}V^{\ast},\quad S^{+}:=Y^{+}\cap
\operatorname*{dom}V^{\ast},\quad S^{-}:=Y^{-}\cap\operatorname*{dom}V^{\ast
},\\
S_{\operatorname{col}}:=Y_{\operatorname{col}}\cap\operatorname*{dom}V^{\ast
},\quad S_{\operatorname{col}}^{+}:=Y_{\operatorname{col}}^{+}\cap
\operatorname*{dom}V^{\ast},\quad S_{\operatorname{col}}^{-}%
:=Y_{\operatorname{col}}^{-}\cap\operatorname*{dom}V^{\ast}.
\end{gather*}
Of course, any of the preceding sets might be empty, $Y_{0}$, $Y^{+}$, $Y^{-}$
being always open, and $Y^{+}$, $Y^{-}$, $Y_{\operatorname{col}}^{+}$,
$Y_{\operatorname{col}}^{-}$ being convex, the convexity of the last two sets
being proved in \cite[Cor.~3]{Zal:18b}. It follows that $S^{+}$, $S^{-}$,
$S_{\operatorname{col}}^{+}$, $S_{\operatorname{col}}^{-}$ are convex, the
first two being open if $\operatorname*{dom}V^{\ast}$ is so; moreover
$\operatorname*{int}S_{\operatorname{col}}^{+}\subset S^{+}$
(resp.$~\operatorname*{int}S_{\operatorname{col}}^{-}\subset S^{-}$) whenever
$\operatorname*{int}S^{+}\neq\emptyset$ (resp.$~\operatorname*{int}S^{-}%
\neq\emptyset$). Obviously,
\begin{align*}
Y^{+}\cup Y^{-}  &  \subset Y_{0}\subset Y_{\operatorname{col}},\quad
Y_{\operatorname{col}}^{+}\cup Y_{\operatorname{col}}^{-}\subset
Y_{\operatorname{col}},\quad Y^{+}=Y_{0}\cap Y_{\operatorname{col}}^{+},\quad
Y^{-}=Y_{0}\cap Y_{\operatorname{col}}^{-},\\
S^{+}\cup S^{-}  &  \subset S_{0}\subset S_{\operatorname{col}},\quad
S_{\operatorname{col}}^{+}\cup S_{\operatorname{col}}^{-}\subset
S_{\operatorname{col}},\quad S^{+}=S_{0}\cap S_{\operatorname{col}}^{+},\quad
S^{-}=S_{0}\cap S_{\operatorname{col}}^{-}.
\end{align*}

In \cite{Zal:18b} we considered a dual function associated to the family
$(q_{k})_{k\in\overline{0,m}}$, which is denoted by $D_{L}$ in this work. More
precisely,%
\[
D_{L}:Y_{\operatorname{col}}\rightarrow\mathbb{R},\quad D_{L}(\sigma
):=L(x,\sigma)\text{ with }A(\sigma)x=b(\sigma).
\]

In a similar way, we consider the (dual objective) function $D$ associated to
$(q_{k})_{k\in\overline{0,m}}$ and $V$ defined by
\[
D:S_{\operatorname{col}}\rightarrow\mathbb{R},\quad D(\sigma):=\Xi
(x,\sigma)\text{ with }A(\sigma)x=b(\sigma);
\]
hence
\begin{equation}
D(\sigma)=D_{L}(\sigma)-V^{\ast}(\sigma)\quad\forall\sigma\in
S_{\operatorname{col}}. \label{r-dld}%
\end{equation}
Setting
\[
x(\sigma):=A(\sigma)^{-1}b(\sigma):=[A(\sigma)]^{-1}\cdot b(\sigma)
\]
for $\sigma\in Y_{0}$, we obtain that
\[
D(\sigma)=\Xi\left(  x(\sigma),\sigma\right)  =-\tfrac{1}{2}\left\langle
b(\sigma),A(\sigma)^{-1}b(\sigma)\right\rangle +c(\sigma)-V^{\ast}%
(\sigma)\quad\forall\sigma\in S_{0}.
\]
From \cite[Prop.~4 (i)]{Zal:18b} we have that $D_{L}$ is concave and upper
semi\-continuous (usc) on $Y_{\operatorname{col}}^{+}$, and convex and lower
semi\-continuous (lsc) on $Y_{\operatorname{col}}^{-}$, and \cite[Eq.~(9)]%
{Zal:18b} holds; moreover $D_{L}(\sigma)$ is attained at any $x\in
\mathbb{R}^{n}$ such that $A(\sigma)x=b(\sigma)$ whenever $\lambda\in
Y_{\operatorname{col}}^{+}\cup Y_{\operatorname{col}}^{-}$, being attained
uniquely at $x:=x(\sigma)$ for $\sigma\in Y^{+}\cup Y^{-}$. Taking into
account (\ref{r-dld}) we have that
\begin{equation}
D(\sigma)=\left\{
\begin{array}
[c]{ccc}%
\min_{x\in\mathbb{R}^{n}}\Xi(x,\sigma) & \text{if} & \sigma\in
S_{\operatorname{col}}^{+},\\
\max_{x\in\mathbb{R}^{n}}\Xi(x,\sigma) & \text{if} & \sigma\in
S_{\operatorname{col}}^{-},
\end{array}
\right.  \label{Pd}%
\end{equation}
the value of $D(\sigma)$ being attained uniquely at $x:=x(\sigma)$ when
$\sigma\in S^{+}\cup S^{-}$ $(\subset S_{0})$; moreover, we have that $D$ is
concave and usc on $S_{\operatorname{col}}^{+}$ as the sum of two concave and
usc functions, while $D$ is a d.c.\ function (difference of convex functions)
on $S_{\operatorname{col}}^{-}$. In general, $D$ is neither convex nor concave
on (the convex set) $S_{\operatorname{col}}^{-}$. Having in view
\cite[Eq.~(11)]{Zal:18b} (or by direct calculations), we have that
\begin{align}
\frac{\partial D}{\partial\sigma_{i}}(\sigma)  &  =\tfrac{1}{2}\left\langle
A(\sigma)^{-1}b(\sigma),A_{i}A(\sigma)^{-1}b(\sigma)\right\rangle
-\left\langle b_{i},A(\sigma)^{-1}b(\sigma)\right\rangle +c_{i}-\frac{\partial
V^{\ast}}{\partial\sigma_{i}}(\sigma)\nonumber\\
&  =\tfrac{1}{2}\left\langle x(\sigma),A_{i}x(\sigma)\right\rangle
-\left\langle b_{i},x(\sigma)\right\rangle +c_{i}-\frac{\partial V^{\ast}%
}{\partial\sigma_{i}}(\sigma)=q_{i}\left(  x(\sigma)\right)  -\frac{\partial
V^{\ast}}{\partial\sigma_{i}}(\sigma) \label{r-dpd}%
\end{align}
for those $\sigma\in\operatorname*{int}S_{0}$ and $i\in\overline{1,m}$ for
which $\frac{\partial V^{\ast}}{\partial\sigma_{i}}(\sigma)$ exists.

\begin{proposition}
\label{p-sub}Assume that $V\in\Gamma(\mathbb{R}^{m})$ is sublinear. Then
$D|_{S_{\operatorname{col}}^{-}}$ is convex; moreover, $\nabla D(\sigma
)=q(x(\sigma))$ for every $\sigma\in S_{0}\cap\operatorname*{int}%
(\operatorname*{dom}V^{\ast})$.
\end{proposition}

Proof. Assume that $V$ is sublinear; it follows that $V^{\ast}=\iota_{\partial
V(0)}$. Then $\Xi(x,\sigma)=L(x,\sigma)\in\mathbb{R}$ for $\sigma
\in\operatorname*{dom}V^{\ast}$, and so $\Xi(x,\cdot)|_{\operatorname*{dom}%
V^{\ast}}$ is convex because $L(x,\cdot)$ is linear and $\operatorname*{dom}%
V^{\ast}$ is a convex set; in particular, $\Xi(x,\cdot
)|_{S_{\operatorname{col}}^{-}}$ is convex because $S_{\operatorname{col}}%
^{-}$ $(\subset\operatorname*{dom}V^{\ast})$ is convex. Using (\ref{Pd}) we
obtain that $D|_{S_{\operatorname{col}}^{-}}$ is convex, too.

Of course, $V^{\ast}$ being constant on $\operatorname*{dom}V^{\ast}$, $\nabla
V^{\ast}(\sigma)=0$ for every $\sigma\in\operatorname*{int}%
(\operatorname*{dom}V^{\ast})$. Taking into account (\ref{r-dpd}), we obtain
that $\nabla D(\sigma)=q(x(\sigma))$ for $\sigma\in S_{0}\cap
\operatorname*{int}(\operatorname*{dom}V^{\ast})$ $(\subset\operatorname*{int}%
S_{\operatorname{col}}^{-})$. \hfill$\square$

\bigskip

Let us denote by $\Gamma_{sc}:=\Gamma_{sc}(\mathbb{R}^{m})$ the class of those
$g\in\Gamma(\mathbb{R}^{m})$ which are essentially strictly convex and
essentially smooth, that is the class of proper lsc convex functions of
Legendre type (see \cite[Section 26]{Roc:72}). Note that any differentiable
and strictly convex function $g:\mathbb{R}^{m}\rightarrow\mathbb{R}$ belongs
to $\Gamma_{sc}(\mathbb{R}^{m});$ moreover, $\Gamma_{sc}(\mathbb{R})$ consists
of those $g\in\Gamma(\mathbb{R})$ which are derivable and strictly convex on
$\operatorname*{int}(\operatorname*{dom}g)$, assumed to be nonempty.

Assume that $g\in\Gamma_{sc}$. Then: $g^{\ast}\in\Gamma_{sc}$,
$\operatorname*{dom}\partial g=\operatorname*{int}(\operatorname*{dom}g)$, and
$g$ is differentiable on $\operatorname*{int}(\operatorname*{dom}g);$
moreover, $\nabla g:\operatorname*{int}(\operatorname*{dom}g)\rightarrow
\operatorname*{int}(\operatorname*{dom}g^{\ast})$ is bijective and continuous
with $\left(  \nabla g\right)  ^{-1}=\nabla g^{\ast}$. In the rest of this
section we assume that $V\in\Gamma_{sc}$, and so $V^{\ast}\in\Gamma_{sc}$,
too. Then, because $V$ is differentiable on $\operatorname*{int}%
(\operatorname*{dom}V)$ and $V^{\ast}$ is differentiable on
$\operatorname*{int}(\operatorname*{dom}V^{\ast})$, clearly
\begin{gather}
\nabla f(x)=A_{0}x-b_{0}+\sum\nolimits_{i=1}^{m}\frac{\partial V}{\partial
y_{i}}\left(  q(x)\right)  \cdot\left(  A_{i}x-b_{i}\right)  \quad\forall x\in
X_{0},\label{r-fd1}\\
\nabla_{\sigma}\Xi(x,\sigma)=q(x)-\nabla V^{\ast}(\sigma)\quad\forall
(x,\sigma)\in\mathbb{R}^{n}\times\operatorname*{int}(\operatorname*{dom}%
V^{\ast}), \label{r-fxids}%
\end{gather}
where
\begin{equation}
X_{0}:=\left\{  x\in\mathbb{R}^{n}\mid q(x)\in\operatorname*{int}%
(\operatorname*{dom}V)\right\}  \subset\operatorname*{dom}f; \label{r-xz}%
\end{equation}
moreover, it follows that (\ref{r-dpd}) holds for $\sigma\in
\operatorname*{int}S_{0}$ $[=S_{0}\cap\operatorname*{int}(\operatorname*{dom}%
V^{\ast})]$ and $i\in\overline{1,m}$, and so%
\begin{equation}
\nabla D(\sigma^{\prime})=q(x(\sigma^{\prime}))-\nabla V^{\ast}(\sigma
^{\prime})=\nabla_{\sigma}\Xi(x(\sigma^{\prime}),\sigma^{\prime})\quad
\forall\sigma^{\prime}\in S_{0}\cap\operatorname*{int}(\operatorname*{dom}%
V^{\ast}). \label{r-fxids3}%
\end{equation}
From (\ref{r-fxids}) and (\ref{r-dpd}) we get
\begin{equation}
\nabla_{\sigma}\Xi(x,\sigma)=0\Longleftrightarrow\big[\sigma\in
\operatorname*{int}S_{0}~\wedge~q(x)=\nabla V^{\ast}(\sigma
)\big ]\Longleftrightarrow\big[x\in X_{0}~\wedge~\sigma=\nabla V\left(
q(x)\right)  \big ]. \label{r-fxids2}%
\end{equation}
From the concavity of $\Xi(x,\cdot)$ for $x\in\mathbb{R}^{n}$ and
(\ref{r-fxids2}) we obtain the next variant of (\ref{r-fb}):
\begin{equation}
f(x)=\sup_{\sigma\in\operatorname*{dom}V^{\ast}}\Xi(x,\sigma)=\sup_{\sigma
\in\operatorname*{int}(\operatorname*{dom}V^{\ast})}\Xi(x,\sigma
)=\Xi\big(x,\nabla V\left(  q(x)\right)  \big)\quad\forall x\in X_{0};
\label{r-fxi2}%
\end{equation}
moreover, using (\ref{r-fd1}) and (\ref{r-fxids}) we obtain that
\begin{equation}
\left[  x\in X_{0}~\wedge~\sigma=\nabla V\left(  q(x)\right)  \right]
\Longrightarrow\left[  \nabla f(x)=\nabla_{x}\Xi(x,\sigma)~\wedge
~f(x)=\Xi(x,\sigma)\right]  . \label{r-nfnxix}%
\end{equation}
Furthermore, using (\ref{r-nxXi}) and (\ref{r-fxids}), for $(x,\sigma
)\in\mathbb{R}^{n}\times\operatorname*{int}(\operatorname*{dom}V^{\ast})$ we
have that
\begin{equation}
\nabla\Xi(x,\sigma)=0\Longleftrightarrow\left[  x\in X_{0}~\wedge
~\sigma=\nabla V\left(  q(x)\right)  ~\wedge~A(\sigma)x=b(\sigma)\right]  .
\label{r17b}%
\end{equation}

\section{The case $\overline{\sigma}\in S_{\operatorname{col}}^{+}$}

The preceding considerations yield directly the next result.

\begin{proposition}
\label{p10}Let $V\in\Gamma(\mathbb{R}^{m})$ and $(\overline{x},\overline
{\sigma})\in\mathbb{R}^{n}\times\operatorname*{dom}V^{\ast}$.

\emph{(i)} Assume that $\nabla_{x}\Xi(\overline{x},\overline{\sigma})=0$ and
$q(\overline{x})\in\partial V^{\ast}(\overline{\sigma})$. Then $(\overline
{x},\overline{\sigma})\in\operatorname*{dom}f\times S_{\operatorname{col}}$,
$\overline{\sigma}\in\partial V\left(  q(\overline{x})\right)  $, and
\begin{equation}
f(\overline{x})=\Xi(\overline{x},\overline{\sigma})=D(\overline{\sigma}).
\label{r-flpd}%
\end{equation}

\emph{(ii)} Moreover, assume that $A(\overline{\sigma})\succeq0$. Then
$\overline{\sigma}\in S_{\operatorname{col}}^{+}$ and
\begin{equation}
f(\overline{x})=\inf_{x\in\operatorname*{dom}f}f(x)=\Xi(\overline{x}%
,\overline{\sigma})=\sup_{\sigma\in S_{\operatorname{col}}^{+}}D(\sigma
)=D(\overline{\sigma}); \label{r-minmax}%
\end{equation}
furthermore, if $\overline{\sigma}\in S^{+}$, then $\overline{x}$ is the
unique global solution of problem $(P)$.
\end{proposition}

Proof. (i) Because $q(\overline{x})\in\partial V^{\ast}(\overline{\sigma})$,
from (\ref{r-dsXi0}) and (\ref{r-fen}) we obtain that
\[
\overline{\sigma}\in\partial V\left(  q(\overline{x})\right)  \quad\wedge\quad
V\left(  q(\overline{x})\right)  +V^{\ast}(\overline{\sigma})=\left\langle
q(\overline{x}),\overline{\sigma}\right\rangle ,
\]
whence $\overline{x}\in\operatorname*{dom}f$ and%
\[
f(\overline{x})=q_{0}(\overline{x})+V\left(  q(\overline{x})\right)
=q_{0}(\overline{x})+\left[  \left\langle q(\overline{x}),\overline{\sigma
}\right\rangle -V^{\ast}(\overline{\sigma})\right]  =\Xi(\overline
{x},\overline{\sigma});
\]
hence the first equality in (\ref{r-flpd}) holds. Because $A(\overline{\sigma
})\overline{x}-b(\overline{\sigma})=\nabla_{x}\Xi(\overline{x},\overline
{\sigma})=0$, we have that $\overline{\sigma}\in S_{\operatorname{col}}$, and
the second equality in (\ref{r-flpd}) holds by the definition of $D$. (ii) By
(i) we have that (\ref{r-flpd}) holds and $\overline{\sigma}\in
S_{\operatorname{col}}$. Because $A(\overline{\sigma})\succeq0$ we have that
$\overline{\sigma}\in S_{\operatorname{col}}^{+}$ and $\Xi(\cdot
,\overline{\sigma})$ is convex, while because $\nabla_{x}\Xi(\overline
{x},\overline{\sigma})=0$ we have that $(f(\overline{x})=)$ $\Xi(\overline
{x},\overline{\sigma})\leq\Xi(x,\overline{\sigma})\leq f(x)$ for
$x\in\operatorname*{dom}f\setminus\{\overline{x}\}$, the latter inequality
being equivalent to
\[
q_{0}(x)+\left\langle q(x),\overline{\sigma}\right\rangle -V^{\ast}%
(\overline{\sigma})\leq q_{0}(x)+V(q(x)),
\]
which is true by the Fenchel--Young inequality [that is the inequality in
(\ref{r-fen})]; furthermore, $\Xi(\overline{x},\overline{\sigma}%
)<\Xi(x,\overline{\sigma})$ when $A(\overline{\sigma})\succ0$. In particular,
$f(\overline{x})=\min_{x\in\operatorname*{dom}f}f(x)$. Using (\ref{r-fb}), the
inclusion $S_{\operatorname{col}}^{+}\subset\operatorname*{dom}V^{\ast}$,
obvious inequalities, and (\ref{Pd}), we get the following sequence of
inequalities:
\begin{align*}
f(\overline{x})  &  =\inf_{x\in\operatorname*{dom}f}f(x)=\inf_{x\in
\operatorname*{dom}f}\sup_{\sigma\in\operatorname*{dom}V^{\ast}}\Xi
(x,\sigma)\geq\inf_{x\in\operatorname*{dom}f}\sup_{\sigma\in
S_{\operatorname{col}}^{+}}\Xi(x,\sigma)\\
&  \geq\sup_{\sigma\in S_{\operatorname{col}}^{+}}\inf_{x\in
\operatorname*{dom}f}\Xi(x,\sigma)=\sup_{\sigma\in S_{\operatorname{col}}^{+}%
}D(\sigma)\geq D(\overline{\sigma}).
\end{align*}
The inequalities above and (\ref{r-flpd}) show that (\ref{r-minmax}) holds.
\hfill$\square$

\begin{proposition}
\label{p10b}Let $V\in\Gamma_{sc}$ and $(\overline{x},\overline{\sigma}%
)\in\mathbb{R}^{n}\times\operatorname*{int}(\operatorname*{dom}V^{\ast})$.

\emph{(i)} Assume that $(\overline{x},\overline{\sigma})$ is a critical point
of $\Xi$. Then $(\overline{x},\overline{\sigma})\in X_{0}\times
S_{\operatorname{col}}$, $\overline{x}$ is a critical point of $f$, and
(\ref{r-flpd}) holds; moreover, if $\overline{\sigma}\in S_{0}$ then
$\overline{\sigma}$ is a critical point of $D$.

\emph{(ii)} Assume that $(\overline{x},\overline{\sigma})$ is a critical point
of $\Xi$ such that $A(\overline{\sigma})\succeq0$. Then $\overline{\sigma}\in
S_{\operatorname{col}}^{+}$ and (\ref{r-minmax}) holds; moreover, if
$A(\overline{\sigma})\succ0$ then $\overline{x}$ is the unique global solution
of problem $(P)$.

\emph{(iii)} Assume that $\overline{\sigma}\in S_{0}$ and $\overline{\sigma}$
is a critical point of $D$. Then $(\overline{x},\overline{\sigma})$ is a
critical point of $\Xi$, where $\overline{x}:=A(\overline{\sigma}%
)^{-1}b(\overline{\sigma})$; therefore, \emph{(i)} and \emph{(ii)} apply.
\end{proposition}

Proof. Observe first that $\partial V^{\ast}(\overline{\sigma})=\{\nabla
V^{\ast}(\overline{\sigma})\}$ because $V^{\ast}$ is differentiable on
$\operatorname*{int}(\operatorname*{dom}V^{\ast})$. (i) Since $\nabla
\Xi(\overline{x},\overline{\sigma})=0$, from (\ref{r-fxids2}) and
(\ref{r-nfnxix}) we have that $q(\overline{x})=\nabla V^{\ast}(\overline
{\sigma})\in\partial V^{\ast}(\overline{\sigma})$, $\overline{x}\in X_{0}$,
and $\nabla f(\overline{x})=\nabla_{x}\Xi(\overline{x},\overline{\sigma})=0$.
Applying Proposition \ref{p10}~(i) we get the first conclusion of (i). Using
(\ref{r-fxids3}) we obtain that $\nabla D(\overline{\sigma})=0$ when
$\overline{\sigma}\in S_{0}$. (ii) As seen in the proof of (i), $q(\overline
{x})\in\partial V^{\ast}(\overline{\sigma})$. The conclusion follows using
Proposition \ref{p10}~(ii).

(iii) Using (\ref{r-fxids3}) we have that $\nabla_{\sigma}\Xi(\overline
{x},\overline{\sigma})=\nabla D(\overline{\sigma})=0$. The choice of
$\overline{x}$ implies $\nabla_{x}\Xi(\overline{x},\overline{\sigma})=0$, and
so $(\overline{x},\overline{\sigma})$ is a critical point of $\Xi$.
\hfill$\square$

\bigskip

In the rest of this section we consider the important particular case in which
$V:=V_{J}:=\iota_{C_{J}}$ for $J\subset\overline{1,m}$, $J^{c}:=\overline
{1,m}\setminus J$, and%
\[
C_{J}:=\{y\in\mathbb{R}^{m}\mid\left[  \forall j\in J:y_{j}=0\right]
~\wedge~\left[  \forall j\in J^{c}:y_{j}\leq0\right]  \}\subset\mathbb{R}%
^{m}.
\]
Of course, $C_{J}$ is a closed convex cone, $V_{J}\in\Gamma(\mathbb{R}^{m})$
is sublinear, and $V_{J}^{\ast}:=\left(  V_{J}\right)  ^{\ast}=\iota
_{\Gamma_{J}}$, where
\[
\Gamma_{J}:=\{\sigma\in\mathbb{R}^{m}\mid\forall j\in J^{c}:\sigma_{j}%
\geq0\};
\]
hence,
\[
\operatorname*{int}\Gamma_{J}:=\{\sigma\in\mathbb{R}^{m}\mid\forall j\in
J^{c}:\sigma_{j}>0\}\neq\emptyset.
\]
For $y,\sigma\in\mathbb{R}^{m}$ we have that
\begin{align}
\sigma\in\partial V_{J}(y)  &  \Longleftrightarrow y\in\partial V_{J}^{\ast
}(\sigma)\Longleftrightarrow\left[  y\in C_{J}~\wedge~\sigma\in\Gamma
_{J}~\wedge~\left\langle y,\sigma\right\rangle =0\right] \label{r-dvdv*}\\
&  \Longleftrightarrow\left[  \forall j\in J:y_{j}=0\right]  ~~\wedge~~\left[
\forall j\in J^{c}:y_{j}\leq0,~\sigma_{j}\geq0,~y_{j}\sigma_{j}=0\right]  .
\label{r-dvdv*b}%
\end{align}
Note that
\[
C_{\overline{1,m}}=\{0\},\quad\Gamma_{\overline{1,m}}=\mathbb{R}^{m},\quad
C_{\emptyset}=\mathbb{R}_{-}^{m}:=-\mathbb{R}_{+}^{m},\quad\Gamma_{\emptyset
}=\mathbb{R}_{+}^{m},
\]
while for $y,\sigma\in\mathbb{R}^{m}$ (\ref{r-dvdv*}) becomes, respectively,%
\begin{gather*}
\sigma\in\partial V_{\overline{1,m}}(y)\Longleftrightarrow y\in\partial
V_{\overline{1,m}}^{\ast}(\sigma)\Longleftrightarrow y=0,\\
\sigma\in\partial V_{\emptyset}(y)\Longleftrightarrow y\in\partial
V_{\emptyset}^{\ast}(\sigma)\Longleftrightarrow\left[  y\in\mathbb{R}_{-}%
^{m}~\wedge~\sigma\in\mathbb{R}_{+}^{m}~\wedge~\left\langle y,\sigma
\right\rangle =0\right]  .
\end{gather*}

For $J\subset\overline{1,m}$ we get $f_{J}:=q_{0}+V_{J}\circ q=q_{0}%
+\iota_{X_{J}}$ and $\Xi_{J}(x,\sigma):=L(x,\sigma)-\iota_{\Gamma_{J}}%
(\sigma)$, where%
\[
X_{J}:=\{x\in\mathbb{R}^{n}\mid\left[  \forall j\in J:q_{j}(x)=0\right]
~\wedge~\left[  \forall j\in J^{c}:q_{j}(x)\leq0\right]
\}=\operatorname*{dom}f_{J}.
\]

So, for $V:=V_{J}$ the problem $(P)$ becomes the problem $(P_{J})$ of
minimizing $q_{0}$ on $X_{J};$ $(P_{\overline{1,m}})$ is the quadratic problem
$(P_{e})$ of minimizing $q_{0}$ on $X_{e}:=X_{\overline{1,m}}$, while
$(P_{\emptyset})$ is the quadratic problem $(P_{i})$ of minimizing $q_{0}$ on
$X_{i}:=X_{\emptyset}$. These problems are considered in \cite{Zal:18b}.

The dual function corresponding to $V_{J}$ is denoted by $D_{J}:=D_{L}%
|_{Y_{\operatorname{col}}^{J}}$, where $Y_{\operatorname{col}}^{J}:=\Gamma
_{J}\cap Y_{\operatorname{col}}$. As observed immediately after getting the
formula of $D$ in (\ref{Pd}), $D_{J}$ is concave on $Y_{\operatorname{col}%
}^{J+}:=\Gamma_{J}\cap Y_{\operatorname{col}}^{+}$, while from Proposition
\ref{p-sub} we have that $D_{J}$ is convex on $Y_{\operatorname{col}}%
^{J-}:=\Gamma_{J}\cap Y_{\operatorname{col}}^{-}$ because $V_{J}$ is sublinear.

\begin{corollary}
\label{cp10}Let $(\overline{x},\overline{\sigma})\in\mathbb{R}^{n}%
\times\mathbb{R}^{m}$ be a $J$-LKKT point of $L$, that is $\nabla
_{x}L(\overline{x},\overline{\sigma})=0$, and
\begin{equation}
\textstyle\left[  \forall j\in J^{c}:\overline{\sigma}_{j}\geq0~\wedge
~\frac{\partial L}{\partial\sigma_{j}}(\overline{x},\overline{\sigma}%
)\leq0~\wedge~\overline{\sigma}_{j}\cdot\frac{\partial L}{\partial\sigma_{j}%
}(\overline{x},\overline{\sigma})=0\right]  ~\wedge~\left[  \forall j\in
J:\frac{\partial L}{\partial\sigma_{j}}(\overline{x},\overline{\sigma
})=0\right]  . \label{r-lkkt}%
\end{equation}
Then $(\overline{x},\overline{\sigma})\in X_{J}\times Y_{\operatorname{col}%
}^{J}$ and $q_{0}(\overline{x})=L(\overline{x},\overline{\sigma}%
)=D_{L}(\overline{\sigma})$. Moreover, assume that $A(\overline{\sigma
})\succeq0$. Then $\overline{\sigma}\in Y_{\operatorname{col}}^{J+}$ and
\[
q_{0}(\overline{x})=\inf_{x\in X_{J}}q_{0}(x)=L(\overline{x},\overline{\sigma
})=\sup_{\sigma\in Y_{\operatorname{col}}^{J+}}D_{L}(\sigma)=D_{L}%
(\overline{\sigma});
\]
furthermore, if $A(\overline{\sigma})\succ0$, then $\overline{x}$ is the
unique global minimizer of $q_{0}$ on $X_{J}$.
\end{corollary}

Proof. From (\ref{r-lkkt}) we have that $\overline{\sigma}\in\Gamma_{J}$, and
so $\nabla_{x}\Xi_{J}(\overline{x},\overline{\sigma})=\nabla_{x}L(\overline
{x},\overline{\sigma})=0$. Since $\frac{\partial L}{\partial\sigma_{j}%
}(\overline{x},\overline{\sigma})=q_{j}(\overline{x})$ for $j\in\overline
{1,m}$, using again (\ref{r-lkkt}) we obtain that $q(\overline{x})\in C_{J}$,
whence $\overline{x}\in X_{J}$, and $\left\langle q(\overline{x}%
),\overline{\sigma}\right\rangle =0$. Using now (\ref{r-dvdv*}) we obtain that
$q(\overline{x})\in\partial V^{\ast}(\overline{\sigma})$. The conclusion
follows now using Proposition \ref{p10} for $V:=V_{J}$. \hfill$\square$

\bigskip

The variant for maximizing $q_{0}$ on $X_{J}$ is the following result; it can
be obtained from the preceding corollary replacing $q_{0}$ by $-q_{0}$ and
$\sigma$ by $-\sigma$ in the definition of $L$.

\begin{corollary}
\label{cp10b}Let $(\overline{x},\overline{\sigma})\in\mathbb{R}^{n}%
\times\mathbb{R}^{m}$ be such that $\nabla_{x}L(\overline{x},\overline{\sigma
})=0$, and
\begin{equation}
\textstyle\left[  \forall j\in J^{c}:\overline{\sigma}_{j}\leq0~\wedge
~\frac{\partial L}{\partial\sigma_{j}}(\overline{x},\overline{\sigma}%
)\leq0~\wedge~\overline{\sigma}_{j}\cdot\frac{\partial L}{\partial\sigma_{j}%
}(\overline{x},\overline{\sigma})=0\right]  ~\wedge~\left[  \forall j\in
J:\frac{\partial L}{\partial\sigma_{j}}(\overline{x},\overline{\sigma
})=0\right]  . \label{l-kktb}%
\end{equation}
Then $(\overline{x},\overline{\sigma})\in X_{J}\times\left(
Y_{\operatorname{col}}\cap(-\Gamma_{J})\right)  $ and $q_{0}(\overline
{x})=L(\overline{x},\overline{\sigma})=D_{L}(\overline{\sigma})$. Moreover,
assume that $A(\overline{\sigma})\preceq0$. Then $\overline{\sigma}\in
Y_{\operatorname{col}}^{J-}:=Y_{\operatorname{col}}^{-}\cap(-\Gamma_{J})$ and
\[
q_{0}(\overline{x})=\sup_{x\in X_{J}}q_{0}(x)=L(\overline{x},\overline{\sigma
})=\inf_{\sigma\in Y_{\operatorname{col}}^{J-}}D_{L}(\sigma)=D_{L}%
(\overline{\sigma}).
\]
furthermore, if $A(\overline{\sigma})\prec0$, then $\overline{x}$ is the
unique global maximizer of $q_{0}$ on $X_{J}$.
\end{corollary}

\begin{corollary}
\label{cp11}Let $(\overline{x},\overline{\sigma})\in\mathbb{R}^{n}%
\times\mathbb{R}^{m}$ be a critical point of $L$. Then $\overline{x}\in X_{e}%
$, $\overline{\sigma}\in Y_{\operatorname{col}}$ and $q_{0}(\overline
{x})=L(\overline{x},\overline{\sigma})=D_{L}(\overline{\sigma})$. Moreover, if
$A(\overline{\sigma})\succeq0$, then $\overline{\sigma}\in
Y_{\operatorname{col}}^{+}$ and
\[
q_{0}(\overline{x})=\inf_{x\in X_{e}}q_{0}(x)=L(\overline{x},\overline{\sigma
})=\sup_{\sigma\in Y_{\operatorname{col}}^{+}}D_{L}(\sigma)=D_{L}%
(\overline{\sigma});
\]
if $A(\overline{\sigma})\preceq0$, then $\overline{\sigma}\in
Y_{\operatorname{col}}^{-}$ and
\[
q_{0}(\overline{x})=\sup_{x\in X_{e}}q_{0}(x)=L(\overline{x},\overline{\sigma
})=\inf_{\sigma\in Y_{\operatorname{col}}^{-}}D_{L}(\sigma)=D_{L}%
(\overline{\sigma}).
\]

\end{corollary}

Proof. For the first two assertions one applies Corollary \ref{cp10} for
$J:=\overline{1,m}$, while for the third assertion one applies Corollary
\ref{cp10b}. \hfill$\square$

\bigskip

Notice that Corollaries \ref{cp10} and \ref{cp10b} are parts of \cite[Prop.~9]%
{Zal:18b} and \cite[Prop.~12]{Zal:18b}, respectively, while Corollary
\ref{cp11} is \cite[Prop.~5\ (i)]{Zal:18b}.

\bigskip

In many papers by DY Gao and his collaborators one speaks about
\textquotedblleft triality theorems\textquotedblright\ in which, besides the
minimax result established for the case $A(\overline{\sigma})\succeq0$ (see
Proposition \ref{p10}), one obtains also \textquotedblleft bi-duality" results
(\textquotedblleft double-min duality" and \textquotedblleft double-max
duality") established for $A(\overline{\sigma})\prec0$, that is $\overline{x}$
and $\overline{\sigma}$ are simultaneously local minimizers (maximizers) for
$f$ on $\operatorname*{dom}f$ and for $D$ on $S^{-}$, respectively.

The next example shows that such triality results are not valid for general
$V\in\Gamma(\mathbb{R}^{m})$, even for $n=m=1$. We concentrate on the case
$\overline{\sigma}\in S^{-}$ of Proposition \ref{p10} (i), that is
$(\overline{x},\overline{\sigma})\in\mathbb{R}^{n}\times\mathbb{R}^{m}$ is
such that $A(\overline{\sigma})\overline{x}=b(\overline{\sigma})$ and
$\overline{\sigma}\in S^{-}\cap\partial V\left(  q(\overline{x})\right)  $,
and so $\overline{x}\in q^{-1}(\operatorname*{dom}V)=\operatorname*{dom}f$.

\begin{example}
\label{ex1}Consider $n:=m:=1$, $V:=\iota_{\mathbb{R}_{-}}$, and $q_{0}%
(x):=-\tfrac{1}{2}x^{2}+x$, $q(x)=q_{1}(x):=\tfrac{1}{2}\left(  x^{2}%
-1\right)  $ for $x\in\mathbb{R}$. Then $f:=f_{\emptyset}=q_{0}+\iota
_{\lbrack-1,1]}$, $A(\sigma)=\sigma-1$, $b(\sigma)=-1$, $c(\sigma)=-\tfrac
{1}{2}\sigma$, whence $L(x,\sigma)=\tfrac{1}{2}(\sigma-1)x^{2}+x-\tfrac{1}%
{2}\sigma$, $Y_{\operatorname{col}}=Y_{0}=\mathbb{R}\setminus\{1\}$,
$x(\sigma)=1/(1-\sigma)$, and so $D_{L}(\sigma)=\tfrac{1}{2}\big(\frac
{1}{1-\sigma}-\sigma\big)$, for $\sigma\in Y_{0};$ moreover, $D:=D_{\emptyset
}=D_{L}$ on $S_{\operatorname{col}}=[0,1)\cup(1,\infty)$. For $\overline
{\sigma}=0$ $[\in S_{\operatorname{col}}^{-}=S^{-}=[0,1)]$ we get
$\overline{x}:=x(0)=1$. Clearly, $0\in\partial V(q(1))=\partial V(0)$
$(=\mathbb{R}_{+})$. Hence the pair $(1,0)$ verifies the hypothesis of
Proposition \ref{p10} (i), even more, $(1,0)$ is a critical point of $L$.
However, by direct verification, or applying Corollary \ref{cp10b}, we obtain
that $\overline{x}=1$ is the unique global maximizer of $f$ on
$\operatorname*{dom}f=[-1,1]$, while applying \cite[Prop.~4 (iv)]{Zal:18b} we
obtain that $\overline{\sigma}=0$ is the unique global minimizer of $D_{L}$ on
$Y_{\operatorname{col}}^{-}$ $[=(-\infty,1)]$, whence $0$ is the unique global
minimizer of $D$ on $S^{-}$. These facts show that \textquotedblleft
double-min duality" and \textquotedblleft double-max duality" are not verified
in the present case.
\end{example}

In DY Gao's works published after 2011 the \textquotedblleft triality
theorems" are established for $V$ a twice differentiable strictly convex
function. Our aim in the sequel is to study the problems of \textquotedblleft
double-min duality" and \textquotedblleft double-max duality" for a special
class of functions $V$. First, in the next section, we establish a result on
positive semidefinite operators in Euclidean spaces needed for getting our
\textquotedblleft bi-duality\textquotedblright\ results.

\section{An auxiliary result}

In order to study the case when $\overline{\sigma}\in S^{-}$, we need the
following result which is probably known, but we have not a reference for it.

\begin{proposition}
\label{p2}Let $X$, $Y$ be nontrivial Euclidean spaces and $H:Y\rightarrow X$
be a linear operator with $H^{\ast}:X\rightarrow Y$ its adjoint. Consider
$Q:=HH^{\ast}:=H\circ H^{\ast}$, $R:=H^{\ast}H$, and
\[
\varphi:X\rightarrow\mathbb{R},~~\varphi(x):=\left\Vert H^{\ast}x\right\Vert
^{2},\quad\psi:Y\rightarrow\mathbb{R},~~\psi(y):=\left\Vert Hy\right\Vert
^{2}.
\]
Then the following assertions hold:

\emph{(a)} $Q$ and $R$ are self-adjoint positive semi-definite operators,
$\ker Q=\ker H^{\ast}$, $\operatorname{Im}Q=\operatorname{Im}H$, $\ker R=\ker
H$, $\operatorname{Im}R=\operatorname{Im}H^{\ast}$; consequently, $H=0$
$\Leftrightarrow$ $Q=0$ $\Leftrightarrow$ $R=0$.

\emph{(b)} Setting $S_{X}:=\{x\in X\mid\left\Vert x\right\Vert =1\}$, one has
$\alpha=\beta$, where
\begin{align}
\alpha &  :=\max_{x\in S_{X}}\varphi(x)=\max\{\lambda\in\mathbb{R}\mid\exists
x\in X\setminus\{0\}:Qx=\lambda x\},\label{r-alfa}\\
\beta &  :=\max_{y\in S_{Y}}\psi(y)=\max\{\lambda\in\mathbb{R}\mid\exists y\in
Y\setminus\{0\}:Ry=\lambda y\}. \label{r-beta}%
\end{align}

\emph{(c)} If $H\neq0$, then $\operatorname{Im}Q\neq\{0\}$, $\operatorname{Im}%
R\neq\{0\}$, and $\gamma=\delta>0$, where
\begin{align}
\gamma &  :=\min_{x\in S_{X}\cap\operatorname{Im}Q}\varphi(x)=\min
\{\lambda>0\mid\exists x\in X\setminus\{0\}:Qx=\lambda x\},\label{r-gama}\\
\delta &  :=\min_{x\in S_{Y}\cap\operatorname{Im}R}\psi(y)=\min\{\lambda
>0\mid\exists y\in Y\setminus\{0\}:Ry=\lambda y\}. \label{r-delta}%
\end{align}

\emph{(d)} The following implications hold:
\begin{gather*}
\min_{x\in S_{X}}\varphi(x)=0\Longleftrightarrow\ker Q\neq
\{0\}\Longleftrightarrow\operatorname{Im}Q\neq X\Longleftrightarrow
\operatorname{Im}H\neq X,\\
\min_{y\in S_{Y}}\psi(y)=0\Longleftrightarrow\ker R\neq
\{0\}\Longleftrightarrow\operatorname{Im}R\neq Y\Longleftrightarrow\ker
H\not =\{0\}.
\end{gather*}

\end{proposition}

Proof. Observe that any result obtained for $Q$ is valid for $R$ because
$(H^{\ast})^{\ast}=H$.

(a) It is obvious that $Q$ is self-adjoint; moreover, $\left\langle
Qx,x\right\rangle =\left\langle HH^{\ast}x,x\right\rangle =\left\langle
H^{\ast}x,H^{\ast}x\right\rangle =\varphi(x)\geq0$ for every $x\in X$. The
inclusions $\ker H^{\ast}\subset\ker Q$ and $\operatorname{Im}Q\subset
\operatorname{Im}H$ are obvious.

Take $x\in\ker Q$, that is $Qx=0$; then $0=\left\langle x,Qx\right\rangle
=\left\Vert H^{\ast}x\right\Vert ^{2}$, and so $x\in\ker H^{\ast}$.

Because $Q$ is self-adjoint, we have that $\operatorname{Im}Q=(\ker Q)^{\perp
}$, and so $X=\ker Q+\operatorname{Im}Q$. Let $x\in\operatorname{Im}H$, that
is $x=Hy$ for some $y\in Y$; then $x=Qu+z$ for some $u\in X$ and $z\in\ker
Q=\ker H^{\ast}$, and so $\left\Vert z\right\Vert ^{2}=\left\langle
z,Hy-HH^{\ast}u\right\rangle =\left\langle H^{\ast}z,y-H^{\ast}u\right\rangle
=0$. It follows that $x=Qu\in\operatorname{Im}Q$.

Because $\ker H=X\Leftrightarrow H=0\Leftrightarrow H^{\ast}=0\Leftrightarrow
\ker H^{\ast}=Y$, the mentioned equivalences follow from the first part.

(b) The conclusion is obvious if $H=0$ (in which case $Q=0$ and $\varphi=0$).
So, let $H\neq0$, and so $Q\neq0$, whence $\alpha>0$. Even if the equalities
in (\ref{r-alfa}) and (\ref{r-beta}) are well known, they will be recovered
below. In fact, the inequalities $\geq$ are almost obvious. Because $\varphi$
is continuous and $S_{X}$ is compact, there exists $\overline{x}\in S_{X}$
such that $\alpha=\varphi(\overline{x})$, and so%
\[
\alpha=\left\Vert H^{\ast}\overline{x}\right\Vert ^{2}=\left\langle
\overline{x},Q\overline{x}\right\rangle \geq\varphi(x)=\left\langle
x,Qx\right\rangle \quad\forall x\in S_{X},
\]
whence $\alpha\left\Vert x\right\Vert ^{2}\geq\left\langle x,Qx\right\rangle
$, or equivalently $\left\langle (\alpha I-Q)x,x\right\rangle \geq0$, for
$x\in X$. Using Schwarz inequality for positive semi-definite operators and
the fact that $\left\langle (\alpha I-Q)\overline{x},\overline{x}\right\rangle
=0$, we get
\[
\left\vert \left\langle (\alpha I-Q)\overline{x},x\right\rangle \right\vert
\leq\sqrt{\left\langle (\alpha I-Q)\overline{x},\overline{x}\right\rangle
}\sqrt{\left\langle (\alpha I-Q)x,x\right\rangle }=0\quad\forall x\in X;
\]
hence $(\alpha I-Q)\overline{x}=0$, that is $Q\overline{x}=\alpha\overline{x}%
$. Hence the inequality $\leq$ holds in (\ref{r-alfa}). Since $Q=HH^{\ast}$,
setting $\overline{y}:=\alpha^{-1/2}H^{\ast}\overline{x}\in Y$, we have that
$\left\Vert \overline{y}\right\Vert =\alpha^{-1/2}\left\Vert H^{\ast}%
\overline{x}\right\Vert =1$, and so $\overline{y}\in S_{Y}$. It follows that
\[
\beta\geq\psi(\overline{y})=\left\Vert H\overline{y}\right\Vert ^{2}%
=\alpha^{-1}\left\Vert HH^{\ast}\overline{x}\right\Vert ^{2}=\alpha
^{-1}\left\Vert Q\overline{x}\right\Vert ^{2}=\alpha^{-1}\left\Vert
\alpha\overline{x}\right\Vert ^{2}=\alpha.
\]
Applying the argument above for $Q$ replaced by $R$, we obtain that
$\alpha\geq\beta$, and so $\alpha=\beta$.

(c) First observe that $S_{X}\cap\operatorname{Im}Q$ is a nonempty compact
set, and so there exists $\overline{x}\in S_{X}\cap\operatorname{Im}Q$ such
that $\gamma=\varphi(\overline{x})$, and so%
\begin{equation}
\gamma=\left\Vert H^{\ast}\overline{x}\right\Vert ^{2}=\left\langle
\overline{x},Q\overline{x}\right\rangle \leq\varphi(x)=\left\langle
x,Qx\right\rangle \quad\forall x\in S_{X}\cap\operatorname{Im}Q. \label{r-gfi}%
\end{equation}
Assuming that $\left\langle \overline{x},Q\overline{x}\right\rangle =0$, as
above, we obtain that $Q\overline{x}=0$, that is $\overline{x}\in\ker Q$.
Since $\ker Q\cap\operatorname{Im}Q=\{0\}$, we get the contradiction $0\in
S_{X}$. Therefore, $\gamma>0$. From (\ref{r-gfi}) we obtain that
$\gamma\left\Vert x\right\Vert ^{2}\leq\left\langle x,Qx\right\rangle $, or
equivalently $\left\langle (Q-\gamma I)x,x\right\rangle \geq0$, for
$x\in\operatorname{Im}Q$. Using Schwarz inequality for the positive
semi-definite operator $\Phi:=(Q-\gamma I)|_{\operatorname{Im}Q}%
:\operatorname{Im}Q\rightarrow\operatorname{Im}Q$ and the fact that
$\left\langle \Phi\overline{x},\overline{x}\right\rangle =0$, we get
$\left\vert \left\langle \Phi\overline{x},x\right\rangle \right\vert \leq
\sqrt{\left\langle \Phi\overline{x},\overline{x}\right\rangle }\sqrt
{\left\langle \Phi x,x\right\rangle }=0$ for all $x\in\operatorname{Im}Q$,
whence $\Phi\overline{x}=0$, that is $Q\overline{x}=\gamma\overline{x}$. As in
the proof of (b) we take $\overline{y}:=\gamma^{-1/2}H^{\ast}\overline{x}%
\in\operatorname{Im}H^{\ast}$; it follows that $\overline{y}\in S_{Y}%
\cap\operatorname{Im}H^{\ast}=S_{Y}\cap\operatorname{Im}R$, and so $\delta
\leq\psi(\overline{y})=\gamma$. The converse inequality follows similarly.

(d) These equivalences are immediate consequences of the equalities in (a) and
the arguments at the beginning of the proof of (c). \hfill$\square$

\section{The case $\overline{\sigma}\in S^{-}$\label{sec-5}}

Throughout this section we assume that $V\in\Gamma_{sc}^{2}$, where
$\Gamma_{sc}^{2}:=\Gamma_{sc}^{2}(\mathbb{R}^{m})$ is the class of those
$g\in\Gamma_{sc}$ which are twice differentiable on $\operatorname*{int}%
(\operatorname*{dom}g)$ with $\nabla^{2}g(y)\succ0$ for $y\in
\operatorname*{int}(\operatorname*{dom}g)$.\footnote{Note that the function
$V$ considered in \cite{VoiZal:13} belongs to $\Gamma_{sc}^{2}(\mathbb{R}).$}
Observe that for $g\in\Gamma_{sc}^{2}$ one has $g^{\ast}\in\Gamma_{sc}^{2}$
and%
\begin{equation}
\nabla^{2}g^{\ast}(\sigma)=\big(\nabla^{2}g\left(  (\nabla g)^{-1}%
(\sigma)\right)  \big)^{-1}\quad\forall\sigma\in\operatorname*{int}%
(\operatorname*{dom}g^{\ast}). \label{r19}%
\end{equation}

In the sequel $V\in\Gamma_{sc}^{2}$. It follows that
\begin{equation}
\left\langle u,\nabla^{2}f(x)u\right\rangle =\bigg\langle u,\big[A_{0}%
+\sum\nolimits_{i=1}^{m}\frac{\partial V}{\partial y_{i}}\left(  q(x)\right)
\cdot A_{i}\big]u\bigg\rangle+\left\langle v_{u},\nabla^{2}V(q(x))v_{u}%
\right\rangle \label{r-fd2}%
\end{equation}
for all $x\in X_{0}$ and $u\in\mathbb{R}^{n}$, where $v_{u}:=\left(
\left\langle u,A_{1}x-b_{1}\right\rangle ,...,\left\langle u,A_{m}%
x-b_{m}\right\rangle \right)  ^{T}$, and
\begin{align*}
\frac{\partial^{2}D}{\partial\sigma_{i}\partial\sigma_{k}}(\sigma)  &
=-\left\langle A_{i}A(\sigma)^{-1}b(\sigma)-b_{i},A^{-1}\left(  A_{k}%
A(\sigma)^{-1}b(\sigma)-b_{k}\right)  \right\rangle -\frac{\partial^{2}%
V^{\ast}}{\partial\sigma_{i}\partial\sigma_{k}}(\sigma)\\
&  =-\left\langle A_{i}x(\sigma)-b_{i},A^{-1}\left(  A_{k}x(\sigma
)-b_{k}\right)  \right\rangle -\frac{\partial^{2}V^{\ast}}{\partial\sigma
_{i}\partial\sigma_{k}}(\sigma)
\end{align*}
for all $\sigma\in\operatorname*{int}S_{0}$ and $i,k\in\overline{1,m}$. It
follows that
\begin{equation}
\big\langle v,\nabla^{2}D(\sigma)v\big\rangle=-\left\langle A_{v}%
x(\sigma)-b_{v},A(\sigma)^{-1}\left(  A_{v}x(\sigma)-b_{v}\right)
\right\rangle -\left\langle v,\nabla^{2}V^{\ast}(\sigma)v\right\rangle
\label{r-d2pd1}%
\end{equation}
for all $v\in\mathbb{R}^{m}$ and $\sigma\in S_{0}$, where
\begin{equation}
A_{v}:=\sum\nolimits_{i=1}^{m}v_{j}A_{j},\quad b_{v}:=\sum\nolimits_{j=1}%
^{m}v_{j}b_{j}\quad(v\in\mathbb{R}^{m}). \label{r-av}%
\end{equation}
The expression above shows that $D$ is strictly concave on $S^{+}$, confirming
the remark done after getting the formulas for $D$ in (\ref{Pd}).

\bigskip

Assume that $(\overline{x},\overline{\sigma})\in X_{0}\times S^{-}$ is a
critical point of $\Xi$; by (\ref{r17b}) we have that $\overline{\sigma
}=\nabla V\left(  q(\overline{x})\right)  $. Because $A(\overline{\sigma
})\prec0$ and $\nabla^{2}V\left(  q(\overline{x})\right)  \succ0$ there exist
non-singular matrices $E\in\mathfrak{M}_{n}$ and $F\in\mathfrak{M}_{m}$ such
that $-A(\overline{\sigma})=E^{\ast}E$ and $\nabla^{2}V\left(  q(\overline
{x})\right)  =F^{\ast}F$, where $E^{\ast}$ and $F^{\ast}$ are the transposed
matrices of $E$ and $F$, respectively; hence $A(\overline{\sigma}%
)^{-1}=-E^{-1}(E^{-1})^{\ast}$ and $\nabla^{2}V^{\ast}\left(  \overline
{\sigma}\right)  =F^{-1}(F^{-1})^{\ast}$. Let us set%
\begin{equation}
d_{i}:=(E^{-1})^{\ast}(A_{i}\overline{x}-b_{i})\in\mathbb{R}^{n}\quad
(i\in\overline{1,m}). \label{r-di}%
\end{equation}
Because for $A_{v}$ defined in (\ref{r-av}) one has
\[
A_{v}\overline{x}-b_{v}=\sum\nolimits_{i=1}^{m}v_{i}(A_{i}\overline{x}%
-b_{i})=\sum\nolimits_{i=1}^{m}v_{i}E^{\ast}d_{i}=E^{\ast}\sum\nolimits_{i=1}%
^{m}v_{i}d_{i},
\]
from (\ref{r-d2pd1}) we obtain that
\begin{align}
\big\langle v,\nabla^{2}D(\overline{\sigma})v\big\rangle  &  =\left\langle
E^{\ast}\sum\nolimits_{i=1}^{m}v_{i}d_{i},E^{-1}(E^{-1})^{\ast}E^{\ast}%
\sum\nolimits_{i=1}^{m}v_{i}d_{i}\right\rangle -\big\langle v,F^{-1}%
(F^{-1})^{\ast}v\big\rangle\nonumber\\
&  =\left\Vert \sum\nolimits_{i=1}^{m}v_{i}d_{i}\right\Vert ^{2}-\left\Vert
(F^{-1})^{\ast}v\right\Vert ^{2}\quad\forall v\in\mathbb{R}^{m}.
\label{r-n2pdef}%
\end{align}
Taking into account (\ref{r-fd2}), we have that
\begin{equation}
\left\langle u,\nabla^{2}f(\overline{x})u\right\rangle =\left\langle
u,A(\overline{\sigma})u\right\rangle +\left\langle v_{u},\nabla^{2}V\left(
q(\overline{x})\right)  v_{u}\right\rangle =\left\Vert Fv_{u}\right\Vert
^{2}-\left\Vert Eu\right\Vert ^{2}\quad\forall u\in\mathbb{R}^{n},
\label{r-n2fef}%
\end{equation}
where%
\begin{equation}
v_{u}=(\left\langle u,A_{i}\overline{x}-b_{i}\right\rangle )_{i\in
\overline{1,m}}=(\left\langle Eu,(E^{-1})^{\ast}(A_{i}\overline{x}%
-b_{i})\right\rangle )_{i\in\overline{1,m}}=(\left\langle Eu,d_{i}%
\right\rangle )_{i\in\overline{1,m}}, \label{r-vu}%
\end{equation}

Let us set
\begin{equation}
J:\mathbb{R}^{m}\rightarrow\mathbb{R}^{n},\quad Jv:=\sum\nolimits_{i=1}%
^{m}v_{i}d_{i}\quad(v\in\mathbb{R}^{m}); \label{r-j}%
\end{equation}
then
\begin{equation}
J^{\ast}:\mathbb{R}^{n}\rightarrow\mathbb{R}^{m},\quad J^{\ast}u=\left(
\left\langle u,d_{1}\right\rangle ,...,\left\langle u,d_{m}\right\rangle
\right)  ^{T}=:\big(\left\langle u,d_{i}\right\rangle _{i\in\overline{1,m}%
}\big)\quad(u\in\mathbb{R}^{n}). \label{r-j*}%
\end{equation}
Take $H:\mathbb{R}^{m}\rightarrow\mathbb{R}^{n}$ defined by $H:=J\circ
F^{\ast}$. Then $H^{\ast}=F\circ J^{\ast}:\mathbb{R}^{n}\rightarrow
\mathbb{R}^{m}$. Because denoting $u^{\prime}:=Eu$ for $u\in\mathbb{R}^{n}$
and $v^{\prime}:=(F^{-1})^{\ast}v$ for $v\in\mathbb{R}^{m}$, from
(\ref{r-n2fef}) and (\ref{r-n2pdef}) we obtain that
\[
\left\langle u,\nabla^{2}f(\overline{x})u\right\rangle =\left\Vert H^{\ast
}u^{\prime}\right\Vert ^{2}-\left\Vert u^{\prime}\right\Vert ^{2}%
,\quad\big\langle v,\nabla^{2}D(\overline{\sigma})v\big\rangle=\left\Vert
Hv^{\prime}\right\Vert ^{2}-\left\Vert v^{\prime}\right\Vert ^{2}.
\]
Because $E$ and $F$ are non-singular, for $\rho\in\{>$, $\geq$, $<$, $\leq\}$
and $\rho^{\prime}\in\{\succ$, $\succeq$, $\prec$, $\preceq\}$ with the
natural correspondence, we have
\begin{gather}
\nabla^{2}f(\overline{x})~\rho^{\prime}~0\Longleftrightarrow\left[  \left\Vert
H^{\ast}u^{\prime}\right\Vert ^{2}~\rho~1~~\forall u^{\prime}\in S_{n}\right]
\Longleftrightarrow\left[  \varphi(u)~\rho~1~~\forall u\in S_{n}\right]
,\label{r-n2f}\\
\nabla^{2}D(\overline{\sigma})~\rho^{\prime}~0\Longleftrightarrow\left[
\left\Vert Hv^{\prime}\right\Vert ^{2}~\rho~1~~\forall v^{\prime}\in
S_{m}\right]  \Longleftrightarrow\left[  \psi(v)~\rho~1~~\forall v\in
S_{m}\right]  , \label{r-n2pd}%
\end{gather}
where $S_{m}:=\{y\in\mathbb{R}^{m}\mid\left\Vert y\right\Vert
=1\}=S_{\mathbb{R}^{m}}$, and $\varphi$, $\psi$ are defined in Proposition
\ref{p2} with
\begin{equation}
H:=J\circ F^{\ast}:\mathbb{R}^{m}\rightarrow\mathbb{R}^{n},\quad H^{\ast
}=F\circ J^{\ast}:\mathbb{R}^{n}\rightarrow\mathbb{R}^{m}. \label{r-H}%
\end{equation}

Recall that $E\in\mathfrak{M}_{n}$ and $F\in\mathfrak{M}_{m}$ are such that
$-A(\overline{\sigma})=E^{\ast}E$ and $\nabla^{2}V\left(  q(\overline
{x})\right)  =F^{\ast}F$, $(d_{i})_{i\in\overline{1,m}}$ are defined in
(\ref{r-di}), $J$ is defined in (\ref{r-j}), and $H$ is defined in (\ref{r-H}).

\bigskip

In the next result we shall use Proposition \ref{p2} for the operator $H$
defined in (\ref{r-H}); therefore, $X=\mathbb{R}^{n}$ and $Y=\mathbb{R}^{m}$.
Using Proposition \ref{p2} if necessary, and setting $\dim\{0\}:=0$, the
following assertions hold:

\begin{itemize}
\item $\dim(\operatorname{Im}H)=\dim(\operatorname{Im}H^{\ast})\leq
\min\{n,m\}$,

\item $\dim(\ker H)+\dim(\operatorname{Im}H)=m$, ~$\dim(\ker H^{\ast}%
)+\dim(\operatorname{Im}H^{\ast})=n$,

\item $\dim(\ker H^{\ast})$ $[=\dim(\ker Q)]$ is equal to the multiplicity of
the eigenvalue $0$ of $Q:=H\circ H^{\ast}$, while $\dim(\ker H)$ is equal to
the multiplicity of the eigenvalue $0$ of $R:=H^{\ast}\circ H$.
\end{itemize}

From the above considerations we obtain the following result.

\begin{proposition}
\label{p-1}Let $(\overline{x},\overline{\sigma})\in X_{0}\times S^{-}$ be a
critical point of $\Xi$. Consider $E\in\mathfrak{M}_{n}$ such that $E^{\ast
}E=-A(\overline{\sigma})$, $d_{i}\in\mathbb{R}^{n}$ $(i\in\overline{1,m})$
defined in (\ref{r-di}), and $H$ defined in (\ref{r-H}).

\emph{(i)} If $\overline{x}$ (resp.\ $\overline{\sigma}$) is a local maximizer
of $f$ (resp.\ $D$), then $\left\Vert Hv\right\Vert \leq1$ for all $v\in
S_{m}$, or, equivalently, $(\alpha=)$ $\beta\leq1$. Conversely, if $\left\Vert
Hv\right\Vert <1$ for all $v\in S_{m}$, then $\overline{x}$ (resp.\ $\overline
{\sigma}$) is a local strict maximizer of $f$ (resp.\ $D$). In particular, if
$A_{i}\overline{x}=b_{i}$ (or equivalently $d_{i}=0$) for all $i\in
\overline{1,m}$, then $\overline{x}$ and $\overline{\sigma}$ are local strict
maximizers of $f$ and $D$, respectively.

\emph{(ii)} If $\overline{x}$ is a local minimizer of $f$, then $\left\Vert
H^{\ast}u\right\Vert \geq1$ for all $u\in S_{n}$; in particular $H$ is
surjective, $m\geq n$, and every positive eigenvalue of $H^{\ast}\circ H$ is
greater than or equal to $1$. Conversely, if $\left\Vert H^{\ast}u\right\Vert
>1$ for all $u\in S_{n}$, then $\overline{x}$ is a local strict minimizer of
$f$; moreover, if $m>n$ then $\overline{\sigma}$ is not a local extremum for
$D$.

\emph{(iii)} If $\overline{\sigma}$ is a local minimizer of $D$, then
$\left\Vert Hv\right\Vert \geq1$ for all $v\in S_{m}$; in particular $H$ is
injective, $m\leq n$, and every positive eigenvalue of $H\circ H^{\ast}$ is
greater than or equal to $1$. Moreover, if $m<n$ then $\overline{x}$ is not a
local extremum for $f$. Conversely, if $\left\Vert Hv\right\Vert >1$ for all
$v\in S_{m}$, then $\overline{\sigma}$ is a local strict minimizer of $D$.

\emph{(iv)} Assume that $m=n$ and $\{A_{i}\overline{x}-b_{i}\mid i\in
\overline{1,m}\}$ is a basis of $\mathbb{R}^{m}$. If $\left\Vert Hv\right\Vert
>1$ for all $v\in S_{m}$, then $\overline{x}$ and $\overline{\sigma}$ are
local strict minimizers of $f$ and $D$, respectively.
\end{proposition}

Proof. Taking into account the well known second order necessary or sufficient
conditions for local extrema of unconstrained problems, the assertions are
immediate consequences of (\ref{r-n2f}), (\ref{r-n2pd}) and Proposition
\ref{p2}. \hfill$\square$

\bigskip

Note that Proposition \ref{p-1}~(iii) gives a positive answer to the question
formulated on the sixth line from below of \cite[p.~234]{VoiZal:13} because in
that case $\overline{\varsigma}$ is a strict local minimum of $P^{d}$ $(=D)$
and $\overline{x}$ is not a local extremum of $P$ $(=f)$ since $m=1<2\leq n$.
In \cite{GaoWu:17} one uses the following assumption: \textquotedblleft(A3)
The critical points of problem ($\mathcal{P}$) are non-singular, i.e., if
$\nabla\Pi(\overline{x})=0$, then $\det\nabla^{2}\Pi(\overline{x})\neq
0$\textquotedblright. Under such a condition we have the following result.

\begin{corollary}
\label{cor-p-1}Let $(\overline{x},\overline{\sigma})\in X_{0}\times S^{-}$ be
a critical point of $\Xi$ such that $\det\nabla^{2}f(\overline{x})\neq0$ [that
is $0$ is not an eigenvalue of $\nabla^{2}f(\overline{x})$]. The following
assertions hold:

\emph{(a)} $\overline{x}$ is a local maximizer of $f$ if and only if
$\left\Vert Hv\right\Vert <1$ for all $v\in S_{m}$, if and only if
$\overline{\sigma}$ is a local maximizer of $D$.

\emph{(b)} Assume that $m=n$. Then $\overline{x}$ is a local minimizer of $f$
if and only if $\left\Vert Hv\right\Vert >1$ for all $v\in S_{m}$, if and only
if $\overline{\sigma}$ is a local minimizer of $D$.
\end{corollary}

Proof. First observe that for $A\in\mathfrak{S}_{n}$ one has $A\succ0$ if and
only if [$A\succeq0\text{ and }\det A\neq0$]. Recall that $\alpha:=\max_{u\in
S_{n}}\left\Vert H^{\ast}u\right\Vert ^{2}=\max_{u\in S_{m}}\left\Vert
Hv\right\Vert ^{2}=:\beta$.

(a) Assume that $\overline{x}$ is a local maximizer of $f$. Then
$A:=\nabla^{2}f(\overline{x})\preceq0$ and so, $A\prec0$. By (\ref{r-n2f}) we
have that $1>\alpha=\beta$ (that is $\left\Vert Hv\right\Vert <1$ for all
$v\in S_{m}$), which at its turn implies that $\overline{\sigma}$ is a local
maximizer of $D$ by Proposition \ref{p-1} (i).

Assume that $\overline{\sigma}$ is a local maximizer of $D$. Then $\alpha
\leq1$ by Proposition \ref{p-1} (i), and so $A\preceq0$ by (\ref{r-n2f}),
whence $A\prec0$. Using again Proposition \ref{p-1} (i), we have that
$\overline{x}$ is a local maximizer of $f$.

The proof of (b) is similar to that of (a). \hfill$\square$

\section{Relations with previous results\label{sec-6}}

In this section we analyze results obtained by DY Gao and his collaborators in
papers dedicated to unconstrained optimization problems, related to
\textquotedblleft triality theorems". The main tool to identify the papers
where this class of problems are considered was to look in the survey papers
\cite{Gao:03b} (which practically includes \cite{Gao:03}), \cite{Gao:09}
(which is almost the same as \cite{Gao:08}), \cite{GaoShe:09} (which is very
similar to \cite{Gao:09b}), \cite{GaoRuaLat:17} (which is the same as
\cite{GaoRuaLat:16}), as well as in the recent book \cite{GaoLatRua:17}.

Though, in order to understand the chronology of the development of this topic
let us quote first the following texts from \cite[p.~40]{GaoRuaLat:17} (see
also \cite[p.~NP30]{GaoRuaLat:16}) and \cite[p.~136]{GaoWu:17} (see also
\cite[p.~5]{GaoWu:12}), respectively:

Q1 -- \textquotedblleft the triality was proposed originally from
post-buckling analysis [42] in \textquotedblleft either-or\textquotedblright%
\ format since the double-max duality is always true but the double-min
duality was proved only in one-dimensional nonconvex analysis
[49]",\footnote{\label{fn-42}The reference \textquotedblleft\lbrack
42]\textquotedblright\ is \textquotedblleft Gao, D.Y.: Dual extremum
principles in finite deformation theory with applications to post-buckling
analysis of extended nonlinear beam theory. Appl.\ Mech.\ Rev.\ 50(11),
S64--S71 (1997)\textquotedblright, while \textquotedblleft\lbrack
49]\textquotedblright\ is DY Gao's book \cite{Gao:00} from our bibliography.
Unfortunately, it is not given the precise place (e.g.\ the page) where the
\textquotedblleft double-min duality\textquotedblright\ was proved for
$n=m=1$.
\par
Among other results, in \textquotedblleft\lbrack42]\textquotedblright\ there
is \textquotedblleft Theorem 7 (Triality Theorem)", stated without proof;
immediately after it is said: \textquotedblleft The proofs of these theorems
are given elsewhere".}

\smallskip Q2 -- \textquotedblleft the triality theorem was formed by these
three pairs of dualities and has been used extensively in nonconvex mechanics
[10, 17] and global optimization [3, 21, 34]. However, it was realized in 2003
[12, 13] that if the dimensions of the primal problem and its canonical dual
are different, the double-min duality (30) needs \textquotedblleft certain
additional conditions\textquotedblright. For the sake of mathematical rigor,
the double-min duality was not included in the triality theory and these
additional constraints were left as an open problem (see Remark 1 in [12],
also Theorem 3 and its Remark in a review article by Gao [13]). By the facts
that the double-max duality (29) is always true and the double-min duality
plays a key role in real-life applications, it was still included in the
triality theory in the either-or form in many applications for the purposes of
perfection in esthesis and some other reasons in reality."\footnote{The
references \textquotedblleft\lbrack3]", \textquotedblleft\lbrack17]" and
\textquotedblleft\lbrack34]" are: \textquotedblleft Fang, S.C., Gao, D.Y.,
Sheu, R.L.,Wu, S.Y.: Canonical dual approach for solving 0--1 quadratic
programming problems. J. Ind.\ Manag.\ Optim.\ 4, 125--142 (2008)",
\textquotedblleft Gao, D.Y., Ogden, R.W.: Multiple solutions to non-convex
variational problems with implications for phase transitions and numerical
computation. Q. J. Mech.\ Appl.\ Math 61(4), 497--522 (2008)", and
\textquotedblleft Ruan, N., Gao, D.Y., Jiao, Y.: Canonical dual least square
method for solving general nonlinear systems of quadratic equations.
Comput.\ Optim.\ Appl.\ 47, 335--347 (2010)", respectively; the references
\textquotedblleft\lbrack10]", \textquotedblleft\lbrack12]", \textquotedblleft%
\lbrack13]", \textquotedblleft\lbrack21]" are our items \cite{Gao:00},
\cite{Gao:03}, \cite{Gao:03b}, and \cite{GaoShe:09}, respectively.}%
~\footnote{On the web-page
http://www.isogop.org/organization/david-y-gao/super-duality-triality
(accessible at least until September 1st, 2014), DY Gao confessed:
\textquotedblleft Actually, I even forgot my this problem left in 2003 [1,2]
due to busy life during those years". So, which is the truth about continuing
to formulate the \textquotedblleft triality theorems" in the \textquotedblleft
either-or form\textquotedblright\ in the papers published in the period
2003--2011?}

\smallskip Having in view Q1 and Q2, it seems that the main steps in the
development of the \textquotedblleft triality theory" are marked by
\cite{Gao:00} (where the triality theorem was proved for the one-dimensional
case), \cite{Gao:03} (where it is mentioned that \textquotedblleft certain
additional conditions\textquotedblright\ are needed for the \textquotedblleft
double-min duality\textquotedblright\ to be valid), \cite{GaoWu:17} and its
preprint version \cite{GaoWu:12a} (where \textquotedblleft this double-min
duality has been proved for ... general global optimization problems", as
mentioned in \cite[p.~40]{GaoRuaLat:17}).

\medskip

Let us compare first our results with those from the most recently published
paper on this topic for general $V$, that is \cite{GaoWu:17}.

\smallskip Putting together Assumptions (A1) and (A2) of Gao and Wu's paper
\cite{GaoWu:17} (see also \cite{GaoWu:12a}), the function $V$ considered there
is real-valued, strictly convex, and twice continuously differentiable on
$\operatorname{Im}q$ (see also \cite[p.~134]{GaoWu:17}). Hence $V$ from
\cite{GaoWu:17} is more general than being in $\Gamma_{sc}^{2}$ when
$\operatorname*{dom}V=\mathbb{R}^{m}.$\footnote{\label{fn1} Note that
$\operatorname*{dom}V=$ $(0,\infty)^{m}$ in \cite[Eq.~(61)]{GaoWu:17}.} Of
course, the strict convexity of $V$ implies $\nabla^{2}V(y)\succeq0$ for
$y\in\mathbb{R}^{m}$, but this property does not imply $(\nabla^{2}V)\left(
q(\overline{x})\right)  \succ0$, which is used for example in \cite[Eq.~(36)]%
{GaoWu:17}. In \textquotedblleft Theorem 2 (Tri-duality Theorem)" (the case
$n=m$) and \textquotedblleft Theorem 3. (Triality Theorem)" (the case
$n\not =m$), $\overline{\sigma}\in S_{\operatorname{col}}$ is a
\textquotedblleft critical point of the canonical problem $(\mathcal{P}^{d})$"
and $\overline{x}:=\left[  A(\overline{\sigma})\right]  ^{-1}b(\overline
{\sigma})$,\footnote{It is not defined what is meant by critical point of the
problem \textquotedblleft$(\mathcal{P}^{d})$ : ext\{$\Pi^{d}(\varsigma
)=-\tfrac{1}{2}\left\langle [G(\varsigma)]^{-1}F(\varsigma),F(\varsigma
)\right\rangle -V^{\ast}(\varsigma)\mid\varsigma\in\mathcal{S}_{a}$\}", where
$\mathcal{S}_{a}$ is our $S_{\operatorname{col}}$ and \textquotedblleft%
$G^{-1}$ should be understood as a generalized inverse if $\det G=0$ [11]",
\textquotedblleft\lbrack11]" being item \cite{Gao:00b} from our bibliography.
Moreover, the formula for $\nabla^{2}\Pi^{d}(\varsigma)$ in Eq.~(34) is not
justified, having in view that $\varsigma\in\mathcal{S}_{a}$
$(=S_{\operatorname{col}}).$} and Assumption (A3) holds, that is $[\nabla
f(x)=0$ $\Rightarrow$ $\det\nabla^{2}f(x)\neq0]$. Our result in the case
$A(\overline{\sigma})\succeq0$ is more general than those in \cite[Ths.~2,
3]{GaoWu:17} not only because the hypothesis on $V$ in Proposition \ref{p10}
is weaker and Assumption (A3) is not present, but also because the conclusion
in \cite[Ths.~2, 3]{GaoWu:17} is weaker, more precisely $f(\overline{x}%
)=\inf_{x\in\mathbb{R}^{n}}f(x)$ $\Leftrightarrow$ $\sup_{\sigma\in
S_{\operatorname{col}}^{+}}D(\sigma)=D(\overline{\sigma})$. In what concerns
the case $A(\overline{\sigma})\prec0$ and $V\in\Gamma_{sc}^{2}$, Corollary
\ref{cor-p-1} is much more precise than the corresponding results in
\cite[Ths.~2, 3]{GaoWu:17} because it is mentioned when $\overline{x}$ and
$\overline{\sigma}$ are local minimizers (maximizers). Moreover, our proofs
are very different from those of \cite{GaoWu:17}, and follow the lines of the
proof of \cite[Prop.~1]{VoiZal:13}.

\medskip

Similar results to those in \cite[Ths.~2, 3]{GaoWu:17} for particular $V$ can
be found in several papers co-authored by DY\ Gao after he became acquainted
with the content of our paper \cite{VoiZal:13}:\footnote{The paper
\cite{VoiZal:13} was submitted to MMOR on 11/09/2009 (manuscript
MMOR-D-09-00165), rejected on DY Gao's report on 15/04/2011, and re-submitted,
without any modification, on 27/04/2011 (manuscript MMOR-D-11-00075); see the
submission date of \cite{GaoWu:12} to arxiv.} \cite{GaoWu:12} (see also
\cite{GaoWu:11}, \cite{GaoWu:11b}), \cite{MorGao:12}, \cite{GaoRuaPar:12},
\cite{RuaGao:14}, \cite{CheGao:16}, \cite{JinGao:16}, \cite{JinGao:17}.

\medskip Gao and Wu in \cite{GaoWu:11}, \cite{GaoWu:11b} and \cite{GaoWu:12}
(which are essentially the same) prove \cite[Ths.~2, 3]{GaoWu:17} for
$V(y):=\tfrac{1}{2}\sum_{k=1}^{m}\beta_{k}y_{k}^{2}$ with $\beta_{k}>0$ and
$b_{k}:=0$ $(k\in\overline{1,m})$ [under Assumption (A3)], using similar
arguments. Note that $\overline{\sigma}$ is taken to be a critical point of
$D$ in \textquotedblleft Theorem 4.3 (Refined Triality
Theorem)\textquotedblright\ (the case $n\neq m$) instead of being a
\textquotedblleft critical point of Problem $(\mathcal{P}^{d})$", as in
\textquotedblleft Theorem 3.1 (Tri-Duality Theorem)".

\medskip Morales-Silva and Gao in \cite{MorGao:11} discuss the problem from
\cite{GaoWu:11} with $A_{0}:=0$ and $m:=1$.
\medskip Morales-Silva and Gao in \cite{MorGao:12} (and \cite{MorGao:15})
consider $V(y):=\sum_{k=1}^{p}\exp(y_{k})+\tfrac{1}{2}\sum_{k=p+1}^{m}%
\beta_{k}y_{k}^{2}$ for $0\leq p\leq m$ (setting $\sum_{k=i}^{j}\gamma_{k}:=0$
when $j<i$) with $\beta_{k}>0$ for $k\in\overline{p+1,m};$ moreover,
$b_{k}:=0$ and $A_{k}\succeq0$ for $k\in\overline{1,m}$ are such that there
exists $(\alpha_{k})_{k\in\overline{1,m}}\subset\mathbb{R}_{+}^{m}$ with
$\sum_{k=1}^{m}\alpha_{k}A_{k}\succ0$. Under Assumption (A3) and using similar
arguments to those in \cite[Th.~2]{GaoWu:12}, they prove \cite[Ths.~2,
3]{GaoWu:17} for $\overline{\sigma}$ \textquotedblleft a stationary point
of\textquotedblright\ $D$.

\medskip Chen and Gao in \cite{CheGao:16} consider $V(y):=\frac{1}{\beta}%
\log\left(  1+\sum_{k=1}^{p}\exp(\beta y_{k})\right)  +\tfrac{1}{2}%
\sum_{k=p+1}^{m}\beta_{k}y_{k}^{2}$ with $\beta,\beta_{k}>0$ $(k\in
\overline{p+1,m})$ and $b_{k}:=0$ $(k\in\overline{1,m})$. By an elementary
computation (and using Cauchy's inequality) one obtains $\nabla^{2}h(z)\succ0$
for $z\in\mathbb{R}^{p}$, where $h(z):=\ln\left(  1+\sum_{k=1}^{p}\exp
(z_{k})\right)  $ for $z\in\mathbb{R}^{p};$ therefore, $V\in\Gamma_{sc}^{2}$.
In \textquotedblleft Theorem 4 (Triality Theorem)\textquotedblright, for
$\overline{\sigma}\in S_{0}$ a critical point of $D$, one obtains the
\textquotedblleft min-max duality" for $\overline{\sigma}\in S^{+}$, while for
$\overline{\sigma}\in S^{-}$ one obtains the \textquotedblleft double-max
duality" and \textquotedblleft double-min duality" (this one for $m=n$)
without using Assumption (A3). However, the proof in the case $\overline
{\sigma}\in S^{-}$ $(=\mathcal{S}_{a}^{-}$ with the notation in
\cite{CheGao:16}) is not convincing.\footnote{\label{fn-chen1}In the proof of
\cite[Th.~4 2.]{CheGao:16} one says: \textquotedblleft Suppose $\overline
{\varsigma}$ is a local maximizer of $\Pi^{d}(\varsigma)$ in $\mathcal{S}%
_{a}^{-}$. Then we have $\nabla^{2}\Pi^{d}(\overline{\varsigma})=-F^{T}%
G_{a}^{-1}F-D^{-1}\preceq0$ and there exists a neighborhood $\mathcal{S}%
_{0}\subset\mathcal{S}_{a}^{-}$ such that for all $\varsigma\in\mathcal{S}%
_{0}$, $\nabla^{2}\Pi^{d}(\varsigma)\preceq0$. Since the map $x=G_{a}^{-1}f$
is continuous over $\mathcal{S}_{a}$, the image of the map over $\mathcal{S}%
_{0}$ is a neighborhood of $\overline{x}$, which we denoted as $\mathcal{X}%
_{0}$. Next we are going to prove that for any $x\in\mathcal{X}_{0}$,
$\nabla^{2}\Pi(x)\preceq0$, which plus the fact that $\overline{x}$ is a
critical point of $\Pi(x)$ implies $\overline{x}$ is a maximizer of $\Pi(x)$
over $\mathcal{X}_{0}$." A similar argument is used for proving \cite[Th.~4
3.]{CheGao:16}, too.
\par
The drawbacks in the quoted text are the following: (a) the fact that
$\nabla^{2}\Pi^{d}(\overline{\varsigma})\preceq0$ implies that $\nabla^{2}%
\Pi^{d}(\varsigma)\preceq0$ for $\varsigma$ in a neighborhood of
$\overline{\varsigma}$ is not motivated ; (b) generally, a continuous function
is not open.} Let us quote Note 1 from \cite[p.~421]{CheGao:16}:
\textquotedblleft We use the same definition of the neighborhood as defined in
[15] (Note 1 on page 306), i.e., a subset $\mathcal{X}_{0}$ is said to be the
neighborhood of the critical point $\overline{x}$ if $\overline{x}$ is the
only critical point in $\mathcal{X}_{0}$."\footnote{\label{fn-chen2}The
reference \textquotedblleft[15]" is item \cite{Gao:03b} from our bibliography.
With this definition, in Example 1 of \cite[p.~426]{CheGao:16}, $\overline
{x}_{3}$ is a minimizer of $\Pi$ on $\mathcal{X}_{0}:=(-1.5,1.5)\setminus
\{\overline{x}_{1},\overline{x}_{2}\};$ of course, this is false.}

\medskip

Jin and Gao in \cite{JinGao:16} and \cite{JinGao:17} (which are essentially
the same) consider practically the same $V$ as in \cite{MorGao:15}, that is
$V(y):=\sum_{k=1}^{p}\exp(y_{k})+\tfrac{1}{2}\sum_{k=p+1}^{m}y_{k}^{2}$ with
$0\leq p\leq m;$ moreover, $b_{k}:=0$ for $k\in\overline{1,m}$, $A_{0}\prec0$
and $A_{k}\succ0$ for $k\in\overline{p+1,m}$. Note that for this $V$, the
statements of \textquotedblleft Theorem 2.\ (Triality theorem)" in
\cite{JinGao:16} and \cite{JinGao:17} and their proofs are almost the same as
those in \cite[Th.~4]{CheGao:16}. The differences are: a) in the case of
\textquotedblleft min-max duality", Assumption 2 in \cite{JinGao:16}
(resp.\ Assumption~1 in \cite{JinGao:17}) implies $\overline{\sigma}\in S^{+}%
$, b) the case $n\neq m$ for the \textquotedblleft double-min duality" is
missing in \cite{JinGao:16} and \cite{JinGao:17}, and c) \textquotedblleft for
some neighborhood $\mathcal{X}_{0}\times\mathcal{S}_{0}\subset\mathbb{R}%
^{n}\times\mathcal{S}_{a}^{-}$ of $(\overline{x},\overline{\zeta})$" from
\cite[Th.~4]{CheGao:16} is replaced by \textquotedblleft$\overline{x}%
\in\mathcal{X}_{0}\subset\mathbb{R}^{n}$ and $\overline{\zeta}\in
\mathcal{S}_{0}\subset\mathcal{S}_{a}^{-}$" in \cite[Th.~2]{JinGao:16} and
\cite[Th.~2]{JinGao:17}.\footnote{Why not take $\mathcal{X}_{0}:=\mathbb{R}%
^{n}$ and $\mathcal{S}_{0}:=\mathcal{S}_{a}^{-}?$} Of course, the drawbacks in
the proof of \cite[Th.~4]{CheGao:16} mentioned in Note \ref{fn-chen1} remain
valid for the proofs of \cite[Th.~2]{JinGao:16} and \cite[Th.~2]%
{JinGao:17}.\footnote{Setting $G:=G(\zeta)$, \cite[Eq.~(38)]{JinGao:16} and
\cite[Eq.~(25)]{JinGao:17} assert that $\eta:=\inf_{x\in\mathbb{R}^{n}}\left(
\tfrac{1}{2}\left\langle x,Gx\right\rangle -\left\langle x,f\right\rangle
\right)  =-\tfrac{1}{2}\left\langle f,G^{-1}f\right\rangle $ if $G\succ0$ and
$\eta=-\infty$ otherwise. In fact $\eta\in\mathbb{R}$ if and only if
$[G\succeq0$~$\wedge~f\in\operatorname{Im}G];$ if $G\succeq0$ and $f=Gx_{0}$
then $\eta=-\tfrac{1}{2}\left\langle x_{0},Gx_{0}\right\rangle $ (see
e.g.~\cite[Prop. 2.1~(i)]{Zal:16}).}

\medskip

A special place among DY Gao's papers published after 2010 is occupied by
\cite{GaoRuaPar:12} and \cite{RuaGao:14}.

\medskip Gao, Ruan and Pardalos in \cite{GaoRuaPar:12} take the same $V$ as in
\cite{GaoWu:12} but Assumption (A3) is not considered. Putting together
Theorems 2 and 3 from \cite{GaoWu:12} for \textquotedblleft$\overline
{\varsigma}$ a critical point of the canonical dual function $P^{d}%
(\overline{\varsigma}),$\textquotedblright\ with the mention \textquotedblleft
If $n\neq m$, the double-min duality (25) holds conditionally", one gets
\textquotedblleft Theorem 2 (Triality Theorem)\textquotedblright\ of
\cite{GaoRuaPar:12}. A detailed proof is provided in the case $\overline
{\varsigma}\in\mathcal{S}_{a}^{+}$ $(=S_{\operatorname{col}}^{+})$. The proof
for the case $\overline{\varsigma}\in\mathcal{S}_{a}^{-}$ $(=S^{-})$ is the following:

\medskip

\textquotedblleft If $\overline{\varsigma}\in\mathcal{S}_{a}^{-}$, the matrix
$G(\overline{\varsigma})$ is a negative definite. In this case, the
Gao--Strang complementary function $\Xi(\overline{x},\overline{\varsigma})$ is
a so-called super-Lagrangian [14], i.e., it is locally concave in both
$x\in\mathcal{X}_{0}\subset\mathcal{X}_{a}$ and $\varsigma\in\mathcal{S}%
_{0}\subset\mathcal{S}_{a}^{-}$. By the fact that

$\max_{x\in\mathcal{X}_{0}}\max_{\varsigma\in\mathcal{S}_{0}}\Xi
(x,\varsigma)=\max_{\varsigma\in\mathcal{S}_{0}}\max_{x\in\mathbb{R}^{n}}%
\Xi(x,\varsigma)\quad(26)$

\noindent holds on the neighborhood $\mathcal{X}_{0}\times\mathcal{S}_{0}$ of
$(\overline{x},\overline{\varsigma})$, we have the double-max duality
statement (24). If $n=m$, we have [33]:

$\min_{x\in\mathcal{X}_{0}}\max_{\varsigma\in\mathcal{S}_{0}}\Xi
(x,\varsigma)=\min_{\varsigma\in\mathcal{S}_{0}}\max_{x\in\mathbb{R}^{n}}%
\Xi(x,\varsigma)\quad(27)$

\noindent which leads to the double-min duality statement (25). This proves
the theorem.\textquotedblright\footnote{The reference \textquotedblleft%
\lbrack14]" is Gao's book \cite{Gao:00}, while \textquotedblleft\lbrack33]" is
\textquotedblleft Gao, D.Y. and Wu, C-Z. (2010). On the Triality Theory in
Global Optimization, to appear in J. Global Optimization (published online
arXiv:1104.2970v1 at http://arxiv.org/abs/1104.2970)"; in fact this paper is
published in another journal (see \cite{GaoWu:12}).}

\bigskip

Ruan and Gao in \cite{RuaGao:14} take the same $V$ as in
\cite{GaoRuaPar:12} [$m$, $(A_{k})$ and $(b_{k})$ being different]
and similarly, Assumption (A3) is not considered. The differences in
\cite[Th.~2]{RuaGao:14} with respect to \cite[Th.~2]{GaoRuaPar:12}
are: (a) $\mathcal{S}_{a}^{+}$ is $S^{+}$ instead of
$S_{\operatorname{col}}^{+}$, (b) \textquotedblleft on the
neighborhood\textquotedblright\ is replaced by \textquotedblleft on
its
neighborhood\textquotedblright, (c) $m=n$ is replaced by $\dim\mathcal{X}%
_{a}=\dim\mathcal{S}_{a}$, and (d) $n\neq m$ in the case
$\overline{\varsigma }\in\mathcal{S}_{a}^{-}$ is missing. In
\cite[Rem.~1]{RuaGao:14} one mentions: \textquotedblleft The
double-max duality statement (24) can be proved easily
by the fact that $\max_{x\in\mathcal{X}_{0}}\max_{\varsigma\in\mathcal{S}_{0}%
}\Xi(x,\varsigma)=\max_{\varsigma\in\mathcal{S}_{0}}\max_{x\in\mathcal{X}_{0}%
}\Xi(x,\varsigma)$ $\forall(x,\varsigma)\in\mathcal{X}_{0}\times
\mathcal{S}_{0}\subset\mathcal{X}_{a}\times\mathcal{S}_{a}^{-}$. The
definition of the neighborhood was given in [32] (Note 2 on p.~479),
i.e.\ $\mathcal{S}_{0}\subset\mathcal{S}_{a}^{-}$ is said to be a
neighborhood of the critical point $\overline{\varsigma}$ if it is
the only critical point of $\Pi^{d}$ on
$\mathcal{S}_{0}$.''\footnote{It is not explained what is meant by
$\dim\mathcal{X}_{a}=\dim\mathcal{S}_{a}$; the reference
\textquotedblleft\lbrack32]\textquotedblright\ is item \cite{Gao:03}
from our bibliography.} (See also our Note \ref{fn-chen2} about this
definition of a neighborhood.)

\medskip

There are very few differences between the proof of \cite[Th.~2]{GaoRuaPar:12}
in the case $\overline{\varsigma}\in\mathcal{S}_{a}^{-}$ and that of
\textquotedblleft Theorem 2 (Triality Theorem)" from \cite{GaoRua:08}, where
$m=1$, for the same case:

\medskip\textquotedblleft If $\overline{\varsigma}\in\mathcal{S}_{a}^{-}$, the
matrix $A_{d}(\overline{\varsigma})$ is negative definite. In this case, the
Gao-Strang complementary function $\Xi(\overline{x},\overline{\varsigma})$ is
a so-called super-Lagrangian (seeGao (2000a)), i.e., it is locally concave in
both $x\in\mathcal{X}_{0}\subset\mathbb{R}^{n}$ and $\varsigma\in
\mathcal{S}_{a}^{-}\subset\mathcal{S}_{a}$. Thus, by the triality theory
developed in Gao (2000a), we have that either

$\min_{x\in\mathcal{X}_{0}}\max_{\varsigma\in\mathcal{S}_{a}^{-}}%
\Xi(x,\varsigma)=\min_{\varsigma\in\mathcal{S}_{0}}\max_{x\in\mathbb{R}^{n}%
}\Xi(x,\varsigma)\quad(23)$

\noindent or

$\max_{x\in\mathcal{X}_{0}}\max_{\varsigma\in\mathcal{S}_{a}^{-}}%
\Xi(x,\varsigma)=\max_{\varsigma\in\mathcal{S}_{0}}\max_{x\in\mathbb{R}^{n}%
}\Xi(x,\varsigma)\quad(24)$

\noindent holds on the neighborhood $\mathcal{X}_{0}\times\mathcal{S}_{0}$ of
$(\overline{x},\overline{\varsigma})$. Thus, the equality (23) leads to the
statement (21), while (24) leads to the statement (22). This proves the
theorem.\textquotedblright

\medskip

It is worth comparing the two proofs above with that (for the same case) of
\textquotedblleft THEOREM 3 (Global Minimizer and Maximizer)\textquotedblright%
\ from \cite{Gao:03} (and of \textquotedblleft Theorem 3 (Global Minimizer and
Maximizer)\textquotedblright\ from \cite{Gao:03b}), where $m=1$:

\medskip\textquotedblleft If $\overline{y}^{\ast}\in\mathcal{Y}_{-}^{\ast}$,
then $(\overline{x},\overline{y}^{\ast})$ is a so-called \emph{super-critical
point} of the extended Lagrangian $\Xi(x,y^{\ast})$, i.e. $\Xi(\overline
{x},\overline{y}^{\ast})$ is locally concave in each of its variables $x$ and
$y^{\ast}$ on the neighborhood $\mathcal{X}_{r}\times\mathcal{Y}_{r}^{\ast}$.
In this case, we have

$P(\overline{x})=\max_{x\in\mathcal{X}_{r}}\max_{y^{\ast}\in\mathcal{Y}%
_{r}^{\ast}}\Xi(x,y^{\ast})=\max_{y^{\ast}\in\mathcal{Y}_{r}^{\ast}}\max
_{x\in\mathcal{X}_{r}}\Xi(x,y^{\ast})=P^{d}(y^{\ast})$

\noindent by the fact that the maxima of the super-Lagrangian $\Xi(x,y^{\ast
})$ can be taken in either order on the open set $\mathcal{X}_{r}%
\times\mathcal{Y}_{r}^{\ast}$ (see [17]). This proves the rest part of the
theorem and (38).\textquotedblright\footnote{The reference \textquotedblleft%
\lbrack17]" is Gao's book \cite{Gao:00}.}

\medskip

The presentation above shows that the papers \cite{GaoRuaPar:12} and
\cite{RuaGao:14} make the transition from the proofs of \textquotedblleft
triality theorems\textquotedblright\ published before 2010 (in which one
observed that for \textquotedblleft$\overline{\varsigma}\in\mathcal{S}_{a}%
^{-}$", $\Xi$ is a \textquotedblleft so-called super-Lagrangian", and so
\textquotedblleft the triality theory developed in Gao (2000a)" applies), and
the proofs of the other ``triality theorems'' published after 2011 with
detailed and complicated (but not completely convincing) proofs for twice
differentiable strictly convex functions $V$.

\medskip

Coming back to Q1, we did not succeed to identify the place\ in \cite{Gao:00}
where \textquotedblleft the double-min duality was proved only in
one-dimensional nonconvex analysis". We may consider the following text from
\cite[p.~131]{GaoWu:17} as a hint for the above assertion:

\medskip Q3 \textquotedblleft Therefore, instead of the mono-duality in static
systems, convex Hamiltonian systems are controlled by the so-called
\emph{bi-duality theory}.

\textbf{Bi-Duality Theorem [10]}: If $(\overline{x},\overline{y}^{\ast})$ is a
critical point of the Lagrangian $L(x,y^{\ast})$, then $\overline{x}$ is a
critical point of $\Pi(x)$, $\overline{y}^{\ast}$ is a critical point of
$\Pi^{\ast}(y^{\ast})$ and $\Pi(\overline{x})=L(\overline{x},\overline
{y}^{\ast})=\Pi^{\ast}(\overline{y}^{\ast})$. Moreover, if $n=m$, we have

$\Pi(\overline{x})=\max_{x\in\mathcal{X}_{k}}\Pi(x)\Leftrightarrow
\max_{y^{\ast}\in\mathcal{Y}_{s}^{\ast}}\Pi^{\ast}(y^{\ast})=\Pi^{\ast
}(\overline{y}^{\ast})\quad(10)$

$\Pi(\overline{x})=\min_{x\in\mathcal{X}_{k}}\Pi(x)\Leftrightarrow
\min_{y^{\ast}\in\mathcal{Y}_{s}^{\ast}}\Pi^{\ast}(y^{\ast})=\Pi^{\ast
}(\overline{y}^{\ast}).\quad(11)$

\noindent This bi-duality is actually a special case of the triality theory in
geometrically linear systems, which was originally presented in Chap.\ 2 [10]
for one-dimensional dynamical systems with a simple proof.\textquotedblright%
\footnote{The reference \textquotedblleft\lbrack10]" is Gao's book
\cite{Gao:00}.}

\medskip

Denoting assertions (10) and (11) above by (65) and (66), respectively, and
putting \textquotedblleft or\textquotedblright\ between them, one gets the
statement of the \textquotedblleft Bi-Duality Theorem\textquotedblright\ from
\cite[p.~148]{GaoWu:17}. Notice that only the \textquotedblleft Bi-Duality
Theorem\textquotedblright\ from \cite[p.~148]{GaoWu:17} is present in the
preprint version of \cite{GaoWu:17}, that is \cite{GaoWu:12a}.\footnote{Notice
the following (easy to be verified) false assertion from
\cite[Acknowledgements]{GaoWu:17}: \textquotedblleft The paper was posted
online on April 15, 2011 at https://arXiv.org/abs/1104.2970"; just compare the
submission dates (and Acknowledgements) of \cite{GaoWu:11} and
\cite{GaoWu:12a}.}

\medskip

The \textquotedblleft Bi-Duality Theory" is presented in \cite[Sect.~2.6.2]%
{Gao:00}. Apparently the result above is related to \textquotedblleft Theorem
2.6.5 (Double-Min and Double-Max Duality)" from \cite[p.~86]{Gao:00} and to
\textquotedblleft Theorem 4 (Bi-Duality Theory [38])\textquotedblright%
\footnote{The reference \textquotedblleft[38]" is Gao's book \cite{Gao:00}.}
from \cite{GaoShe:09}; in these two theorems there are not references to the
dimensions. Example 4.5 from \cite{StrVoiZal:11} provides a counterexample for
both \cite[Th.~2.6.5]{Gao:00} and \cite[Th.~4]{GaoShe:09}, as well as for the
bi-duality theorem from \cite[p.~131]{GaoWu:17}; however, that example is not
a counterexample for the bi-duality theorem from \cite[p.~148]{GaoWu:17}.

\medskip Another hint should be \cite[Sect.~3.5]{Gao:00} which is called
\textquotedblleft Tri-Extremum Principles and Triality
Theory\textquotedblright, with its subsections 3.5.2 and 3.5.3 which are
called \textquotedblleft Triality Theorems\textquotedblright\ and
\textquotedblleft Tri-Duality Theory", respectively. At the beginning of
\cite[Sect.~3.5]{Gao:00} it is said:

\smallskip\textquotedblleft In this section we present the so-called triality
theory under the following assumption. Assumption 3.5.1 Let $\{(\mathcal{U}%
,\mathcal{U}^{\ast});\left\langle \ast,\ast\right\rangle \}$ and
$\{(\mathcal{E},\mathcal{T});\left\langle \ast,\ast\right\rangle \}$ be two
inner product spaces. ...

\noindent(A1) $\Lambda:I\times\mathcal{U}\rightarrow\mathcal{E}$ is a
quadratic operator $\Lambda(u)=\tfrac{1}{2}a(x)u^{\prime}(x)^{2}%
+b(x)u^{\prime}(x)+c(x)$, $a(x)>0$ $\forall x\in I$, where $a,b,c\in
\mathcal{C}^{1}(I)$ are given real-valued functions;

\noindent(A2) $F:\mathcal{U}_{a}\subset\mathcal{U}\rightarrow\mathbb{R}$ is a
linear, G\^{a}teaux differentiable functional and, on $\mathcal{U}_{a}%
\times\mathcal{U}_{a}^{\ast}\subset\mathcal{U}\times\mathcal{U}^{\ast}$,
$u^{\ast}=DF(u)\Leftrightarrow u=DF^{c}(u^{\ast})\Leftrightarrow\left\langle
u,u^{\ast}\right\rangle =F(u)+F^{c}(u^{\ast})$;

\noindent(A3) $\bar{W}:\mathcal{E}_{a}\subset\mathcal{E}\rightarrow\mathbb{R}$
is either convex or concave and on $\mathcal{E}_{a}\times\mathcal{T}%
_{a}\subset\mathcal{E}\times\mathcal{T}$, the Legendre duality relations
$\varsigma=D\bar{W}(\xi)\Leftrightarrow\xi=D\bar{W}^{c}(\varsigma
)\Leftrightarrow\left\langle \xi;\varsigma\right\rangle =\bar{W}(\xi)+\bar
{W}^{c}(\varsigma)$ hold."

\smallskip From \cite[(3.107)]{Gao:00}, \cite[(3.108)]{Gao:00} and
\cite[(3.113)]{Gao:00} we learn that $\Pi(u)=\bar{W}(\Lambda(u))-F(u)$ for
$u\in\mathcal{U}_{k}$ with $\mathcal{U}_{k}=\{u\in\mathcal{U}_{a}\mid
\Lambda(u)\in\mathcal{E}_{a}\}$, $L(u,\varsigma)=\left\langle \Lambda
(u);\varsigma\right\rangle )-\bar{W}^{c}(\varsigma)-F(u)$, and $\Pi
^{d}(\varsigma)=F^{c}(u^{\ast}(\sigma))-\bar{W}^{c}(\varsigma)-G^{c}%
(\varsigma)$, respectively, in which $F^{c}(u^{\ast}(\sigma))$ is the Legendre
conjugate of $F(u)$, and $G^{c}:\mathcal{T}_{\varnothing}\rightarrow
\mathbb{R}$ is a pure complementary gap functional.

\smallskip

The above text shows that, at least in \cite[Sect.~3.5]{Gao:00}, $\mathcal{U}$
is a function space like $H^{1}(I)$. Of course, $F$ being a linear function on
$\mathcal{U}_{a}$ $(\subset\mathcal{U})$, $\mathcal{U}_{a}$ has to be a linear
subspace endowed we the trace topology. A linear functional $f$ defined on a
topological vector space $U$ is G\^{a}teaux differentiable if and only if $f$
is continuous, in which case $Df(u)=f$ for every $u\in U;$ moreover, it is not
possible to speak about \textquotedblleft the Legendre conjugate of $F$". So,
(A2) has not a mathematical meaning. Moreover, in order to speak about
$D\bar{W}(\xi)$ and $D\bar{W}^{c}(\varsigma)$ in (A3), one needs
$\mathcal{E}_{a}$ and $\mathcal{T}_{a}$ be at least algebraically open
(convex) subsets of $\mathcal{E}$ and $\mathcal{T}\mathbf{,}$ respectively. It
is clear that the concerned spaces are not one-dimensional.

Because \textquotedblleft Theorem 3.5.2 (Triality Theorem)" from \cite{Gao:00}
does not refer to primal and dual functions as in the usual formulations of
\textquotedblleft triality theorems" we quote such a result from \cite{Gao:09}
(which is maybe the last one) attributed to (Gao, 2000a), that is our
reference \cite{Gao:00}.

\medskip

\textquotedblleft\emph{Theorem 3} (Triality theory (Gao, 2000a)). Suppose that
$\overline{\varsigma}$ is a critical point of $P^{d}$ and $\overline
{x}=G^{\dag}(\overline{\varsigma})\tau(\overline{\varsigma})$. If
$G(\overline{\varsigma})\succeq0$, then $\overline{x}$ is a global minimizer
of $(\mathcal{P})$, $\overline{\varsigma}$ is a global maximizer of
$(\mathcal{P}^{d})$, and $\min_{x\in\mathcal{X}_{_{a}}}P(x)=\Xi(\overline
{x},\overline{\varsigma})=\max_{\varsigma\in\mathcal{S}_{c}^{+}}%
P^{d}(\varsigma)$. If $G(\overline{\varsigma})\prec0$, then on a neighborhood
$\mathcal{X}_{o}\times\ \mathcal{S}_{o}\subset\mathcal{X}_{a}\times
\mathcal{S}_{c}^{-}$ of $(\overline{x},\overline{\varsigma})$, we have either
$\min_{x\in\mathcal{X}_{o}}P(x)=\Xi(\overline{x},\overline{\varsigma}%
)=\min_{\varsigma\in\mathcal{S}_{o}}P^{d}(\varsigma)$, or $\max_{x\in
\mathcal{X}_{o}}P(x)=\Xi(\overline{x},\overline{\varsigma})=\min_{\varsigma
\in\mathcal{S}_{o}}P^{d}(\varsigma).$"

\medskip We consider that there is a misprint in the last $\min_{\varsigma
\in\mathcal{S}_{o}}P^{d}(\varsigma)$ of \cite[Th.~3]{Gao:09}; it has to be
replaced by $\max_{\varsigma\in\mathcal{S}_{o}}P^{d}(\varsigma)$, as in
\cite[Th.~2]{Gao:08} (and all the other Gao's papers containing a
\textquotedblleft triality theorem").

In \cite[Th.~3]{Gao:09}, \textquotedblleft$\mathcal{X}_{_{a}}\subset
\mathbb{R}^{n}$ is a given feasible space", and \textquotedblleft without
losing much generality", $V:\mathcal{E}_{a}\rightarrow\mathbb{R}$
\textquotedblleft is convex and lower semi\-continuous\textquotedblright.
Moreover \textquotedblleft$G^{\dag}$ is the Moore--Penrose generalized inverse
of $G$". Without looking to details, \cite[Th.~3]{Gao:09} is similar to
\textquotedblleft Theorem 3.5.3 (Tri-Duality Theorem)" from \cite{Gao:00};
note that the Moore--Penrose generalized inverse is not considered in
\cite{Gao:00}.

It is worth quoting the most recent version of the general \textquotedblleft
triality theorem", that is \cite[Th.~3]{GaoRuaLat:17} (the same as
\cite[Th.~3]{GaoRuaLat:16}):

\medskip\textquotedblleft Theorem 3 (Triality theorem) Suppose $\overline{\xi
}^{\ast}$ is a stationary point of $\Pi^{d}(\xi^{\ast})$ and $\overline{\chi
}=G(\overline{\xi}^{\ast})^{-1}\overline{\xi}^{\ast}$. If $\overline{\xi
}^{\ast}\in S_{c}^{+}$, we have

$\Pi(\overline{\chi})=\min_{\chi\in\mathcal{X}_{c}}\Pi(\chi)\Leftrightarrow
\max_{\xi^{\ast}\in\mathcal{S}_{c}^{+}}\Pi^{d}(\xi^{\ast})=\Pi^{d}%
(\overline{\xi}^{\ast});\quad(30)$

\noindent If $\overline{\xi}^{\ast}\in S_{c}^{-}$, then on a neighborhood$^{5}%
$ $\mathcal{X}_{o}\times\mathcal{S}_{o}\subset\mathcal{X}_{c}\times
\mathcal{S}_{c}^{-}$ of $(\overline{\chi},\overline{\xi}^{\ast})$, we have either

$\Pi(\overline{\chi})=\max_{\chi\in\mathcal{X}_{o}}\Pi(\chi)\Leftrightarrow
\max_{\xi^{\ast}\in\mathcal{S}_{o}}\Pi^{d}(\xi^{\ast})=\Pi^{d}(\overline{\xi
}^{\ast}),\quad(31)$

\noindent or (only if $\dim\overline{\chi}=\dim\overline{\xi}^{\ast}$)

$\Pi(\overline{\chi})=\min_{\chi\in\mathcal{X}_{o}}\Pi(\chi)\Leftrightarrow
\min_{\xi^{\ast}\in\mathcal{S}_{o}}\Pi^{d}(\xi^{\ast})=\Pi^{d}(\overline{\xi
}^{\ast}).\quad(32)$\textquotedblright

\bigskip

Note 5 in \cite[Th.~3]{GaoRuaLat:17} (quoted above) is \textquotedblleft The
neighborhood $\mathcal{X}_{o}$ of $\overline{\chi}$ means that on which,
$\overline{\chi}$ is the only stationary point" (see also our Note
\ref{fn-chen2}). Related to this theorem, in \cite[p.~14, 15]{GaoRuaLat:17}
(and \cite[p.~NP13]{GaoRuaLat:16}) it is said:

\medskip\textquotedblleft The triality theory was first discovered by Gao 1996
in post-buckling analysis of a large deformed beam [42, 52]. The
generalization to global optimization was made in 2000 [51]. It was realized
in 2003 that the double-min duality (32) holds under certain additional
condition [57, 58]. Recently, it is proved that this additional condition is
simply $\dim\overline{\chi}=\dim\overline{\xi}^{\ast}$ to have the strong
canonical double-min duality (32), otherwise, this double-min duality holds
weakly in subspaces of $\mathcal{X}_{o}\times\mathcal{S}_{o}$ [79, 80, 112,
113].\textquotedblright\footnote{Compare this text with Q2. The references
\textquotedblleft\lbrack52]\textquotedblright\ and \textquotedblleft%
\lbrack112]\textquotedblright\ are \textquotedblleft Gao, D.Y.: Finite
deformation beam models and triality theory in dynamical post-buckling
analysis. Int.\ J. Non-Linear Mech.\ 5, 103--131 (2000)" and \textquotedblleft
Morales-Silva, D.M., Gao, D.Y.: Complete solutions and triality theory to a
nonconvex optimization problem with double-well potential in $\mathbb{R}^{n}$.
Numer.\ Algebra Contr.\ Optim.\ 3(2), 271--282 (2013)\textquotedblright, for
\textquotedblleft\lbrack42]\textquotedblright\ see Note \ref{fn-42}, while
\textquotedblleft\lbrack51]\textquotedblright, \textquotedblleft%
\lbrack57]\textquotedblright\ , \textquotedblleft\lbrack58]\textquotedblright,
\textquotedblleft\lbrack79]\textquotedblright, \textquotedblleft%
\lbrack80]\textquotedblright\ and \textquotedblleft\lbrack
113]\textquotedblright\ are the items \cite{Gao:00b}, \cite{Gao:03b},
\cite{Gao:03}, \cite{GaoWu:12}, \cite{GaoWu:17} and \cite{MorGao:15} from our
bibliography, respectively. Reference \textquotedblleft\lbrack
112]\textquotedblright\ seems to be the published version of \cite{MorGao:11}%
.}

\bigskip Coming back to \cite[Th.~3]{GaoRuaLat:17}, we have to know which are
the conditions on the function corresponding to our $V$, that is
$\Phi$. At the beginning of Section \textquotedblleft2.4 Triality
Theory\textquotedblright\ of \cite{GaoRuaLat:17} it is said
\textquotedblleft
we need to assume that the canonical function $\Phi:\mathcal{E}_{a}%
\rightarrow\mathbb{R}$ is convex". In \cite[Def.~2]{GaoRuaLat:17} it
is said: \textquotedblleft A real-valued function
$\Phi:\mathcal{E}_{a}\rightarrow \mathbb{R}$ is called canonical if
the duality mapping $\partial
\Phi:\mathcal{E}_{a}\rightarrow\mathcal{E}_{a}^{\ast}$ is one-to-one
and onto\textquotedblright, while on \cite[p.~10]{GaoRuaLat:17} it
is said: \textquotedblleft A canonical function $\Phi(\xi)$ can also
be nonsmooth but should be convex such that its conjugate can be
well-defined by Fenchel transformation
$\Phi^{\sharp}(\xi^{\ast})=\sup\{\left\langle \xi,\xi^{\ast
}\right\rangle
-\Phi(\xi)\mid\xi\in\mathcal{E}_{a}\}$.\textquotedblright\ This
means that $V:=\Phi\in\Gamma(\mathbb{R}^{m})$, that is $V$ is the
same as in \cite[Th.~3]{Gao:09}. However, the hypotheses of
\cite[Th.~3]{GaoRuaLat:17} are stronger than those of
\cite[Th.~3]{Gao:09} because in the latter one asks
$\overline{\sigma}\in S_{0}$ (instead of $S_{\operatorname{col}}$)
and, for having the \textquotedblleft double-min
duality\textquotedblright, one assumes that
$\dim\overline{\chi}=\dim\overline{\xi}^{\ast}$.\footnote{Which is
the meaning of
\textquotedblleft$\dim\overline{\chi}=\dim\overline{\xi}^{\ast}$"?
Why is not \cite[Th.~3]{GaoRuaLat:17} attributed to \cite{Gao:00} at
least for $n=m=1?$}

\medskip So, the framework of \cite[Th.~3]{Gao:09} is that of Proposition
\ref{p10}; however, applying the latter we obtain only the first
assertions of \cite[Th.~3]{Gao:09} and \cite[Th.~3]{GaoRuaLat:17}.
Example 19 from \cite{Zal:18b} shows that the \textquotedblleft
double-max duality\textquotedblright\ and \textquotedblleft
double-min duality\textquotedblright\ are not true for $n=2$ and
$m\in\{1,2\}$, taking $V:=\iota_{\{0\}}\in\Gamma(\mathbb{R})$ for
$m:=1$ and $V:=\iota
_{\mathbb{R}_{-}^{2}}\in\Gamma(\mathbb{R}_{-}^{2})$ for $m:=2$.
Moreover, for $V$ from Example \ref{ex1}, in which $D$ $(=P^{d})$ is
differentiable on its domain $[0,1)\cup(1,\infty)\}$, one has that
$D^{\prime}(0)=0$ and $0\succ G(0)$ $[=A(0)]$, but
$\overline{x}:=x(0)$ and $(\overline{\sigma}=)$
$\overline{\varsigma}:=0$ are not simultaneously local minimizers
(maximizers) for $P$ $(=f)$ and $P^{d}$ on $\mathcal{X}_{_{a}}$
$(=[-1,1])$ and $[0,1)$ $(=\mathcal{S}_{c}^{-})$, respectively. In
particular, Example \ref{ex1} shows that the assertion
\textquotedblleft double-max duality is always true" from Q1 is
false. In particular, even Theorems 3 in \cite{Gao:03} and
\cite{Gao:03b} are false because \textquotedblleft double-max
duality\textquotedblright\ is false, as mentioned above.

\medskip Because the proofs for the \textquotedblleft bi-duality" given after
2011, less that in \cite{GaoRuaPar:12}, are sufficiently involved and refer to
$V$ in a restricted class of convex functions, the natural question is what is
happening with the \textquotedblleft bi-duality\textquotedblright\ results
when $V\in\Gamma_{sc}(\mathbb{R}^{m})$. So, we formulate the following open problem:

\medskip

\textbf{Open problem.} Is the next statement true? Let $V\in\Gamma
_{sc}(\mathbb{R}^{m})$, $\overline{\sigma}\in S^{-}\cap\operatorname*{int}%
(\operatorname*{dom}V^{\ast})$ be a critical point of $D$, and $\overline
{x}:=A(\overline{\sigma})^{-1}b(\overline{\sigma})$. Then $\overline{x}$ is a
local maximizer of $f$ on $\operatorname*{dom}f$ if and only if $\overline
{\sigma}$ is a local maximizer of $D$ on $S^{-};$ moreover, if $m=n$ then
$\overline{x}$ is a local minimizer of $f$ on $\operatorname*{dom}f$ if and
only if $\overline{\sigma}$ is a local minimizer of $D$ on $S^{-}%
$.\footnote{Recall Q2 where it is said: \textquotedblleft these additional
constraints were left as an open problem (see Remark 1 in [12], also Theorem 3
and its Remark in a review article by Gao [13])". In fact in \cite{Gao:03} and
\cite{Gao:03b} there are not open problems related to CDT. In Mathematical
Economics there is an interesting axiom denoted NFL, and coming from
\textquotedblleft no free lunch"; this could be translated by `one gets
nothing from nothing', so for getting even the \textquotedblleft double-max
duality\textquotedblright\ one needs \textquotedblleft certain additional
conditions".} \medskip

In this context let us quote from \cite[p.~40]{GaoRuaLat:17} (or
\cite[NP~30]{GaoRuaLat:16}) and \cite[p.~19]{Gao:16}, respectively:

\smallskip

Q4 \textquotedblleft Six papers are in this group on the triality theory. By
listing simple counterexamples (cf.\ e.g., [137]), Voisei and Zalinescu
claimed: \textquotedblleft a correction of this theory is impossible without
falling into trivial\textquotedblright.$^{11}$ However, even some of these
counterexamples are correct, they are not new. This type of counterexamples
was first discovered by Gao in 2003 [57, 58], i.e., the double-min duality
holds under certain additional constraints (see Remark on page 288 [57] and
Remark 1 on page 481 [58]). But neither [57] nor [58] was cited by Voisei and
Zalinescu in their papers. ... $^{11}$This sentence is deleted by Voisei and
Zalinescu in their revision of [137] after they were informed by referees that
their counterexamples are not new and the triality theory has been
proved."\footnote{The references \textquotedblleft\lbrack137]",
\textquotedblleft\lbrack57]" and \textquotedblleft\lbrack58]" are
\cite{VoiZal:11}, \cite{Gao:03} and \cite{Gao:03b} from our bibliography,
respectively.} \footnote{Quite detailed answers to this kind of assertions can
be found in \cite[Sect.~2]{Zal:16b}.} \footnote{From the text
\textquotedblleft even some of these counterexamples are correct, they are not
new" we have to understand that all our counterexamples are not new, and only
some of them are correct. So, the authors of \cite{GaoRuaLat:17} had to
mention explicitly those counterexamples which are not correct.}

\smallskip

Q5 \textquotedblleft Regarding the so-called \textquotedblleft not convincing
proof\textquotedblright, serious researcher should provide either a convincing
proof or a disproof, rather than a complaint\textquotedblright.

\medskip

Paraphrasing the text in Q5, we could say: \emph{Regarding the text in Q4, as
a serious and honest researcher, DY Gao should have mentioned either his
results which are not true, or even more, he should have written down those
\textquotedblleft additional constraints\textquotedblright\ under which the
conclusions of those results become true, rather than the complaint that
Voisei and Zalinescu never cited either \cite{Gao:03} or \cite{Gao:03b}}.

\section{Conclusions}

-- In Proposition \ref{p10} we showed, with a simple proof, that the
\textquotedblleft min-max duality" from the \textquotedblleft
triality theorem\textquotedblright\ for problem $(P)$ is true for
$V$ a proper lower semi\-continuous convex function on
$\mathbb{R}^{m}$. Moreover, we showed that the \textquotedblleft
min-max duality" from quadratic minimization problems with quadratic
constraints can be obtained using Proposition \ref{p10}.

-- We pointed out which are the relationships between the facts that
$\overline{\sigma}$ and $\overline{x}:=x(\overline{\sigma})$ are
(strict) local maximizers (minimizers) of $D$ on $S^{-}$ and of $f$
on $\operatorname*{dom}f$, respectively, in the case in which
$V\in\Gamma _{sc}^{2}$ (see Proposition \ref{p-1}). In particular,
in Corollary \ref{cor-p-1}, we recovered Theorems 2 and 3 of
\cite{GaoWu:17} under \cite[Assumption 3]{GaoWu:17} for
$V\in\Gamma_{sc}^{2}$, our result being less precise that of
\cite[Th.~3]{GaoWu:17} for $m\neq n;$ however, see the discussion
about \cite{GaoWu:17} in Section \ref{sec-6}.

-- In Section \ref{sec-6} we compared our results with those on
\textquotedblleft triality theorems\textquotedblright\ published by DY Gao and
his collaborators after 2010, mentioning several drawbacks in proofs and
inconsistencies in statements and presentations.

-- We showed that the \textquotedblleft double-min duality\textquotedblright%
\ and \textquotedblleft double-max duality\textquotedblright\ of the
general \textquotedblleft triality theorem\textquotedblright\ from
\cite{GaoRuaLat:17} are false even for $m=n=1$. -- We formulated an
open problem concerning the \textquotedblleft double-min
duality\textquotedblright\ and \textquotedblleft double-max
duality\textquotedblright\ when $V\in\Gamma_{sc}$, problem related
to that mentioned in Q2.

\medskip

\textbf{Acknowledgement} We thank Prof.\ Marius Durea for reading a previous
version of the paper and for his useful remarks.

\end{document}